\documentclass[a4paper,12pt]{article}
\title{Full analysis of the Green's function for a
singularly perturbed
convection-diffusion
       problem 
       \\in three dimensions\footnote{
       This work has been supported by Science Foundation Ireland
       under the Research Frontiers Programme 2008;
       Grant 08/RFP/MTH1536}}
\author{Sebastian~Franz\footnote{
         Institut f\"ur Numerische Mathematik,
         Technische Universit\"at Dresden,
         01062 Dresden, Germany\newline
         e-mail: sebastian.franz@tu-dresden.de}
        \and
        Natalia~Kopteva\footnote{
         Department of Mathematics and Statistics,
         University of Limerick,
         Limerick,
         Ireland\newline
         e-mail: natalia.kopteva@ul.ie}
        }

\usepackage{manfnt}

\usepackage{amsmath}
\usepackage{amsfonts}
\usepackage{amssymb}
\usepackage{amsthm}

\usepackage{graphicx}

\newcommand{\ts}{\textstyle}
\newcommand{\ds}{\displaystyle}
\newcommand{\pt}{\partial}

\newcommand{\E}{\mathrm{e}}

\newcommand{\cube}{{\mbox{\scriptsize\mancube}}}

\newcommand{\laplace}{\Delta}
\newcommand{\grad}{\nabla}
\newcommand{\eps}{\varepsilon}
\newcommand{\R}{\mathbb{R}}
\newcommand{\norm}[2]{\|{#1}\|_{#2}}
\newcommand{\bignorm}[2]{\left\| {#1} \right\|  _{#2}}

\renewcommand{\tilde}{\widetilde}
\newcommand{\ve}[1]{\boldsymbol{#1}}

\numberwithin{equation}{section}

\theoremstyle{plain}
\newtheorem{thm}{Theorem}[section]
\newtheorem{lem}[thm]{Lemma}

\newtheorem{cor}[thm]{Corollary}
\newtheorem{rem}[thm]{Remark}

\begin{document}

  \maketitle
   \begin{abstract}
A linear singularly perturbed convection-diffusion
  problem with characteristic layers is
considered in three dimensions.
Sharp bounds for the associated Green's function and its derivatives
are established in the $L_1$ norm.
The dependence of these bounds on the small perturbation parameter
is shown explicitly.
The obtained estimates will be used in a forthcoming numerical analysis of the considered
problem.

The present article is a more detailed version of our recent paper
\cite{FK09_2}.



   \end{abstract}

   \section{Introduction}
%
   In this paper we consider
   the following problem posed in the unit-cube domain
   $\Omega=(0,1)^3$:%
   \begin{subequations}\label{eq:Lu}
   \begin{align}
     \label{eq:Lu_a}
     L_{\ve{x}}u(\ve{x})
       =-\eps\,\laplace_{\ve{x}}u(\ve{x})-\pt_{x_1}(a(\ve{x})\,u(\ve{x}))+b(\ve{x})\,u(\ve{x})
      &=f(\ve{x})
     \quad \mbox{for }\ve{x}\in\Omega,\\
     u(\ve{x})&=0\,\hspace*{1.0cm}\mbox{for }\ve{x}\in\partial\Omega.
   \end{align}
   \end{subequations}
   Here $\eps$ is a small positive parameter, and we assume that the coefficients
   $a$ and $b$  are sufficiently smooth ($a,\,b\in
   C^\infty(\bar\Omega)$).
   We also assume,  for some positive constant $\alpha$, that
   \begin{gather}\label{assmns}
      a(\ve{x})\geq \alpha>0,
      \qquad b(\ve{x})-\pt_{x_1}a(\ve{x})\geq 0
      \qquad\mbox{for~all~}\ve{x}\in\bar\Omega.
   \end{gather}
   Under these assumptions, \eqref{eq:Lu_a} is a singularly
   perturbed elliptic equation, also referred to as a
   convection-dominated convection-diffusion equation.
Its solutions typically exhibits sharp interior and boundary layers.
   This equation
   serves as a model for Navier-Stokes equations at large Reynolds
   numbers or (in the linearised case) of Oseen equations and
   provides an excellent paradigm for numerical techniques in the
   computational fluid dynamics \cite{RST08}.

\begin{figure}[tb]
      \centerline{
      \includegraphics[width=0.45\textwidth]{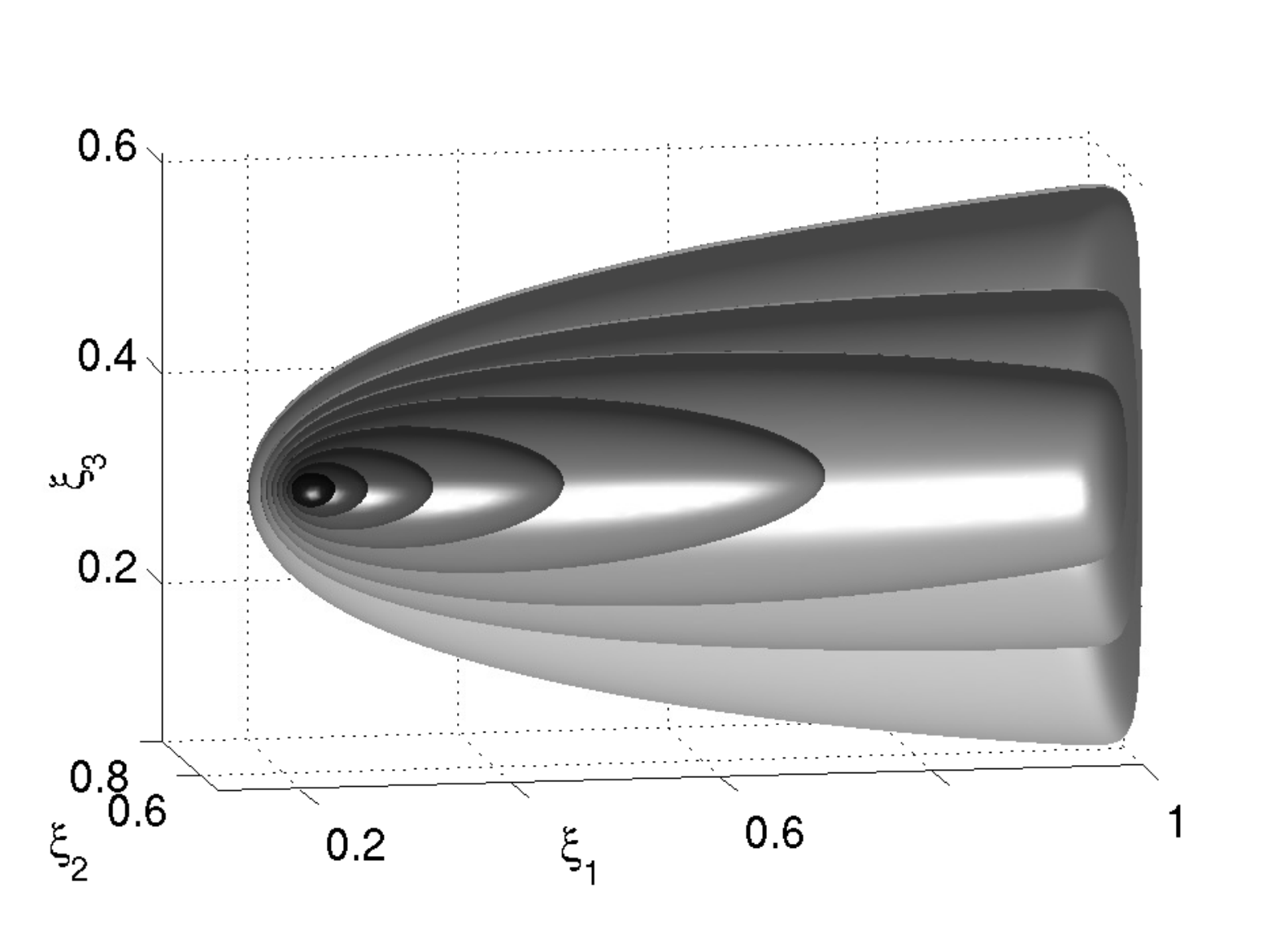}
      \includegraphics[width=0.45\textwidth]{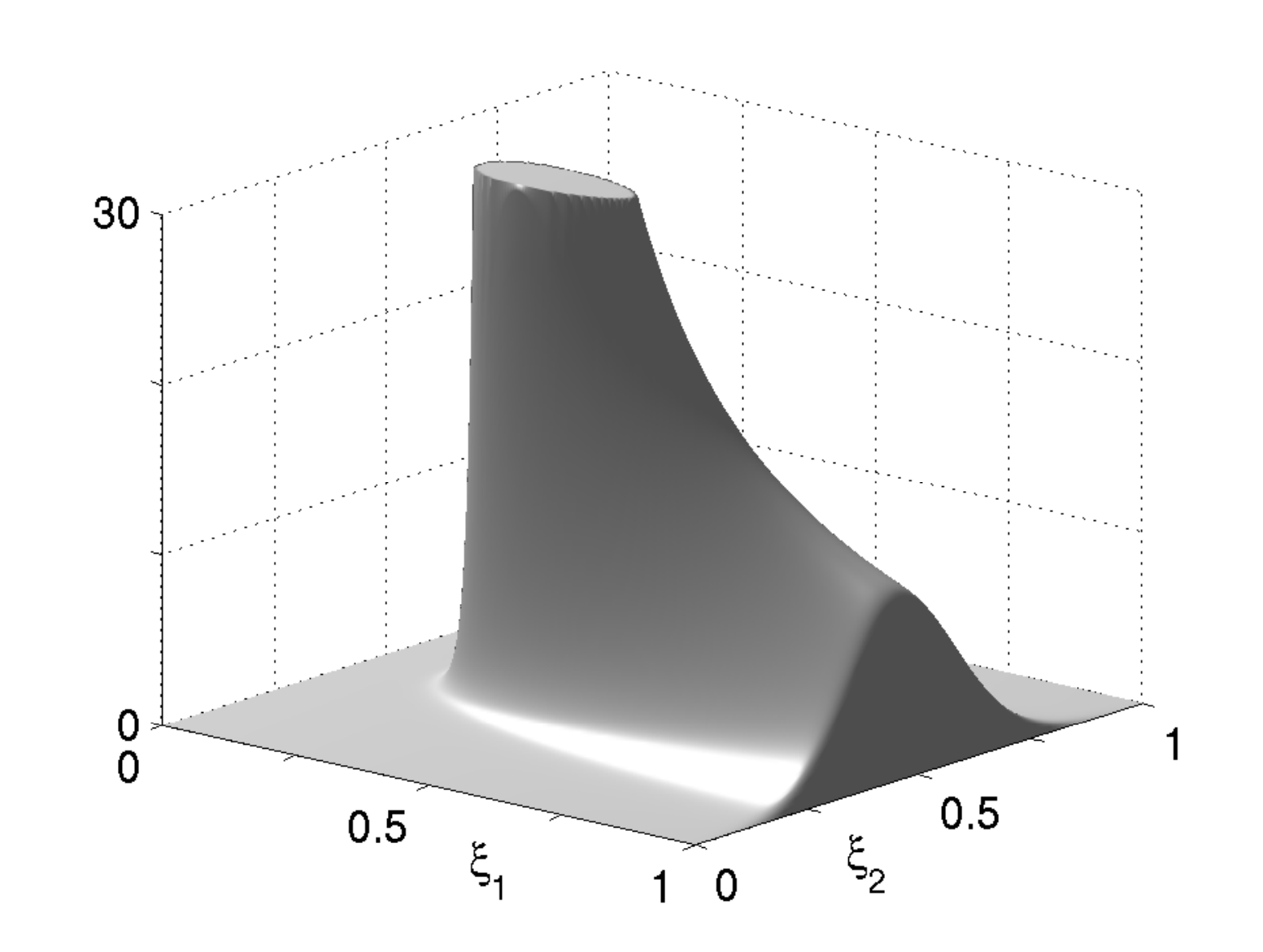}}
      \caption{Anisotropy of the Green's function $G$ associated with
      \eqref{eq:Lu} for $\eps=0.01$ and $\ve{x}=(\frac15,\frac12,\frac13)$.
      Left: 
      isosurfaces at values of
      $1,\,4,\,8,\,16,\,32,\,64,\,128$,  and $256$.
      Right: a two-dimensional graph for fixed $\xi_3=x_3$.}
      \label{fig:green}
   \end{figure}

   The Green's function for the convection-diffusion problem~\eqref{eq:Lu}
   exhibits a strong anisotropic structure, which is demonstrated by Figure~\ref{fig:green}.
   This reflects the complexity of solutions
   of this problem; it should be noted that
   problems of this type
   require an intricate asymptotic analysis
   \cite[Section IV.1]{Ilin}, \cite{KSh87}; see also \cite[Chapter~IV]{Shi92},
   \cite[Chapter~III.1]{RST08} and \cite{KSt05,KSt07}.
   We also refer the reader to
   D\"orfler \cite{Dorf99}, who, for a similar problem,
   gives extensive a priori solution estimates.

   Our interest in considering the Green's function of problem
   \eqref{eq:Lu} and estimating its derivatives is motivated by the
   numerical analysis of this computationally challenging problem.
   More specifically, we shall use the obtained estimates in the
   forthcoming paper \cite{FK10_NA} to derive robust a posteriori
   error bounds for computed solutions of this problem  using
   finite-difference methods. (This approach is related to recent
   articles \cite{Kopt08,CK09}, which address the numerical solution
   of singularly perturbed equations of reaction-diffusion type.)
   In a more general numerical-analysis context, 
   we note that sharp estimates for continuous Green's functions (or
   their generalised versions) frequently play a crucial role in a
   priori and a posteriori error analyses \cite{erikss,Leyk,notch}.

   The purpose of the present paper is to establish
   sharp bounds for the derivatives of the Green's function  in the
   $L_1$ norm (as they will be used to estimate the error in the
   computed solution in the dual $L_\infty$ norm \cite{FK10_NA}).
   Our estimates will be  \emph{uniform in the small perturbation
   parameter $\eps$} in the sense that any dependence on $\eps$ will
   be shown explicitly.
   Note also that our estimates will be {\it sharp}
   (in the sense of Theorem~\ref{theo_lower})
   up to an $\eps$-independent constant multiplier.
   We employ the analysis technique used in
   \cite{FK10_1},
which we now extend to a three-dimensional problem. Roughly
speaking, we freeze the coefficients and estimate the
corresponding explicit frozen-coefficient Green's function, and
then we investigate the difference between the original and the
frozen-coefficient Green's functions.
This procedure is often called the parametrix
method.
To make this paper more readable, we deliberately follow some of
the notation and presentation of \cite{FK10_1}.

   The paper is organised as follows.
   In Section~\ref{sec:def}, the Green's
   function associated with problem \eqref{eq:Lu} is defined and
   upper bounds for its derivatives are stated in Theorem~\ref{thm:main},
   the main result of the paper.
   The corresponding lower bounds are then given in Theorem~\ref{theo_lower}.
   In Section~\ref{sec:def_green_const}, we obtain the fundamental solution
  for a constant-coefficient version of \eqref{eq:Lu} in the domain $\Omega=\R^3$.
  This fundamental solution is bounded in
   Section~\ref{sec:bounds_green_general}.
   It is then used in Section~\ref{sec:bounds_green_const}
   to construct certain
   approximations of the frozen-coefficient Green's functions
   for the domains $\Omega=(0,1)\times\R^2$ and $\Omega=(0,1)^3$.
   The difference between these approximations  and the original
   variable-coefficient Green's function
   is estimated in Section~\ref{sec:main_proof}, which completes the proof of
   Theorem~\ref{thm:main}.
\medskip

   \noindent\textit{Notation.} Throughout the paper, $C$, as well as $c$, denotes
   a generic positive constant
   that may take different values in different formulas,
   but is {\it independent of the small diffusion coefficient $\eps$}.
   A subscripted $C$ (e.g., $C_1$) denotes a positive constant
   that takes a fixed value, and is also independent of $\eps$.
   The usual Sobolev spaces $W^{m,p}(D)$ and $L_p(D)$ on any measurable
   domain $D\subset\R^3$ are used. The $L_p(D)$-norm is denoted by $\norm{\cdot}{p\,;D}$
   while the $W^{m,p}(D)$-norm is denoted by $\norm{\cdot}{m,p\,;D}$.
   By $\ve{x}=(x_1,x_2,x_3)$ we denote an element in $\R^3$.
   For an open ball centred at $\ve{x'}$ of radius $\rho$, we use the notation
   $B(\ve{x'},\rho)=\{\ve{x}\in\R^3: \sum_{k=1,2,3}(x_k-x'_k)^2<\rho^2\}$.
   The notation $\pt_{x_m} f$, $\pt^2_{x_m} f$ and $\laplace_{\ve{x}}$
   is employed for the
   first- and second-order partial derivatives
   of a function $f$ in variable $x_m$, and
 the Laplacian in variable $\ve{x}$, respectively, while
 $\pt^2_{x_k x_m} f$ will denote a mixed derivative of $f$.

%
\section{Definition of Green's function. Main result}\label{sec:def}
%
   The Green's function $G=G(\ve{x};\ve\xi)$
   associated with \eqref{eq:Lu}, satisfies, for each fixed $\ve{x}\in\Omega$,
   \begin{subequations}\label{eq:Green_adj}
   \begin{align}
     \hspace{-0.3cm}
     L^*_{\ve\xi}G(\ve{x};\ve\xi)
          :=-\eps\,\laplace_{\ve\xi}G
            +a(\ve\xi)\,\pt_{\xi_1}G
            +b(\ve\xi)\,G
         &=\delta(\ve{x}-\ve\xi)
     &&\mbox{for}\;\ve\xi\in\Omega,\\
     G(\ve{x};\ve\xi)&=0&&\mbox{for}\;\ve\xi\in\partial\Omega.
   \end{align}
   \end{subequations}
   Here $L^*_{\ve\xi}$ is the adjoint differential operator to $L_{\ve{x}}$,
   and $\delta(\cdot)$ is the three-dimensional Dirac $\delta$-distribution.
   The unique solution $u$ of \eqref{eq:Lu} allows the representation
   \begin{gather}\label{eq:sol_prim}
     u(\ve{x})=\iiint_{\Omega}G(\ve{x};\ve\xi)\,f(\ve\xi)\,d\ve\xi\,.
   \end{gather}
   It should be noted that the Green's function $G$ also satisfies, for each fixed $\ve\xi\in\Omega$,
   \begin{subequations}\label{eq:Green_prim}
   \begin{align}
     L_{\ve{x}}G(\ve{x};\ve\xi)
         =-\eps\,\laplace_{\ve{x}}G
          -\pt_{x_1}(a(\ve{x})\,G)
          +b(\ve{x})\,G
        &=\delta(\ve{x}-\ve\xi)&&\mbox{for}\; \ve{x}\in\Omega,\\
     G(\ve{x};\ve\xi)&=0&&\mbox{for}\;\ve{x}\in\partial\Omega.
   \end{align}
   \end{subequations}
   Consequently, the unique solution $v$ of the adjoint problem
   \begin{subequations}\label{eq:Lu_adj}
   \begin{align}
     L^*_{\ve{x}}v(\ve{x})=-\eps\,\laplace_{\ve{x}}v
     +b(\ve{x})\,\pt_{x_1}v+c(\ve{x})\,v&=f(\ve{x})
     \quad \mbox{for}\;\ve{x}\in\Omega,\\
     v(\ve{x})&=0\hspace*{1.0cm}\,\mbox{for}\;\ve{x}\in\partial\Omega,
   \end{align}
   \end{subequations}
   is given by
   \begin{gather}\label{eq:sol_adj}
     v(\ve\xi)=\iiint_{\Omega}G(\ve{x};\ve\xi)\,f(\ve{x})\,d\ve{x}\,.
   \end{gather}

   We start with a preliminary result for $G$.

   \begin{lem}\label{lem_G_L1}
    Under assumptions \eqref{assmns}, the Green's function $G$
    associated with \eqref{eq:Lu} satisfies
     \begin{gather}\label{G_L1}
       \iint_{(0,1)^2}\! |G(\ve{x};\ve{\xi})|\,d\xi_2\, d\xi_3
         \leq C
       ,\qquad\quad
       \|G(\ve{x};\cdot)\|_{1\,;\Omega}
         \leq C,
     \end{gather}
     where $C$ is some positive $\eps$-independent constant.
   \end{lem}
   \begin{proof}
      The first estimate of \eqref{G_L1} is given in the proof of
      \cite[Theorem 2.10]{Dorf99} (see also \cite[Theorem~III.1.22]{RST08}
      and \cite{And_2003} for similar two-dimensional results).
      The second desired estimate follows.
   \end{proof}
   We now state the main result of this paper.

   \begin{thm}\label{thm:main}
      The Green's function $G$ associated with \eqref{eq:Lu},\,\eqref{assmns}
      in the unit-cube domain $\Omega=(0,1)^3$ satisfies, for all $\ve{x}\in\Omega$,
      the following bounds
      \begin{subequations}
      \begin{align}
         \norm{\pt_{\xi_1} G(\ve{x};\cdot)}{1;\Omega}
            &\leq C(1+|\ln \eps|),\label{eq:thm:G_xi}\\
         \norm{\pt_{\xi_k} G(\ve{x};\cdot)}{1;\Omega}+
         \norm{\pt_{x_k} G(\ve{x};\cdot)}{1;\Omega}
            &\leq C\eps^{-1/2},\quad k=2,\,3.\label{eq:thm:G_eta}
      \end{align}
      Furthermore, for any ball $B(\ve{x}',\rho)$ of radius $\rho$ centered at
      any $\ve{x}'\in\bar\Omega$, we have
      \begin{align}
         \norm{G(\ve{x};\cdot)}{1,1;B(\ve{x}',\rho)}
            &\leq C\eps^{-1}\rho,\label{eq:thm:G_grad}
      \end{align}
      while for the ball $B(\ve{x},\rho)$ of radius $\rho$ centered at $\ve{x}$
      we have
      \begin{align}
         \norm{\pt^2_{\xi_1} G(\ve{x};\cdot)}{1;\Omega\setminus B(\ve{x},\rho)}
            &\leq C\eps^{-1}\ln(2+\eps/\rho),\label{eq:thm:G_xixi}\\
         \norm{\pt^2_{\xi_k} G(\ve{x};\cdot)}{1;\Omega\setminus B(\ve{x},\rho)}
            &\leq C\eps^{-1}(|\ln\eps|+\ln(2+\eps/\rho)),
            \quad k=2,\,3.\label{eq:thm:G_etaeta}
      \end{align}
      \end{subequations}
      Here $C$ is some positive $\eps$-independent constant.
   \end{thm}
We devote the rest of the paper to the proof of this theorem, which will be completed in
Section~\ref{sec:main_proof}.

 In view of the solution representation \eqref{eq:sol_prim},
 Theorem~\ref{thm:main} yields
 a number of a priori solution estimates for our original problem.
E.g., the bounds
   \eqref{eq:thm:G_xi}, \eqref{eq:thm:G_eta} immediately imply the
   following result.
   \begin{cor}\label{cor_apriori}
      Let $f(\ve{x})=\pt_{x_1} F_1(\ve{x})+\pt_{x_2} F_2(\ve{x})+\pt_{x_3} F_3(\ve{x})$
      with $F_1,\,F_2,\,F_3\in L_\infty(\Omega)$.
      Then there exists a unique solution $u\in L_\infty(\Omega)$
      of problem \eqref{eq:Lu}, \eqref{assmns}, for which
      we have the bound
      \begin{gather}\label{apriori}
         \norm{u}{\infty\,;\Omega}
            \leq C\bigl[\,(1+|\ln\eps|)\,\norm{F_1}{\infty\,;\Omega}
                        + \eps^{-1/2}\,(\norm{F_2}{\infty\,;\Omega}+
                                        \norm{F_3}{\infty\,;\Omega})\,\bigr].
      \end{gather}
   \end{cor}

   It can be anticipated from an
   inspection of the bounds for an explicit fundamental solution in a
   constant-coefficient case (see Section~\ref{sec:bounds_green_general})
   that the upper estimates of Theorem~\ref{thm:main} are {\it sharp}.
   Indeed, one can prove the following result.

   \begin{thm}[\cite{FK09_lower}]\label{theo_lower}
      Let $\eps\in(0,c_0]$ for some sufficiently small positive $c_0$.
      Set $a(\ve{x}):=\alpha$ and $b(\ve{x}):=0$ in \eqref{eq:Lu}.
      Then the Green's function $G$ associated with this problem
     in the unit cube $\Omega=(0,1)^3$ satisfies,
      for all $\ve{x}\in[\frac14,\frac34]^3$,
      the following lower bounds:
      \begin{subequations}\label{eq_thm:main_lower}
      \begin{align}
         \norm{\pt_{\xi_1} G(\ve{x};\cdot)}{1;\Omega}
            &\geq c\,|\ln \eps|,\\
         \norm{\pt_{\xi_k} G(\ve{x};\cdot)}{1;\Omega}
            &\geq c\,\eps^{-1/2},\,k=2,\,3.
       \intertext{Furthermore, for any ball $B(\ve{x};\rho)$ of radius $\rho\le\frac18$, we have}
         \norm{G(\ve{x};\cdot)}{1,1;\Omega\cap B(\ve{x};\rho)}
            &\geq \begin{cases}
                    c\,\rho/\eps, & \mbox{for~}\rho\le 2\eps,\\
                    c\,(\rho/\eps)^{1/2},&\mbox{otherwise},\\
                  \end{cases}\\
         \norm{\pt^2_{\xi_1} G(\ve{x};\cdot)}{1;\Omega\setminus B(\ve{x};\rho)}
                       &\geq c\,\eps^{-1}\ln(2+\eps/\rho),
           \quad&&\mbox{for~}\rho\le c_1\eps,
            \\
            \norm{\pt^2_{\xi_k} G(\ve{x};\cdot)}{1;\Omega\setminus B(\ve{x};\rho)}
            &\geq c\,\eps^{-1}(\ln(2+\eps/\rho)+|\ln\eps|)
            &&\mbox{for~}\rho\le{\textstyle\frac18},\,k=2,\,3.
      \end{align}
      Here $c$ and $c_1$ are $\eps$-independent positive constants.
     \end{subequations}
   \end{thm}
%

%
   \section{Fundamental solution in the constant-coefficient case}\label{sec:def_green_const}
%
In our analysis, we invoke the observation that constant-coefficient versions of
the two problems \eqref{eq:Green_adj} and \eqref{eq:Green_prim}
   that we have for $G$, can be easily solved explicitly
when posed in $\R^3$.
   So in this section we shall explicitly solve simplifications
   of  \eqref{eq:Green_adj} and \eqref{eq:Green_prim}. To get these simplifications,
   we employ the parametrix method and so
   freeze the coefficients in these problems by
   replacing $a(\ve{\xi})$ by $a(\ve{x})$ in \eqref{eq:Green_adj}, and
   replacing $a(\ve{x})$ by $a(\ve{\xi})$ in \eqref{eq:Green_prim},
   and also setting $b:=0$; the frozen-coefficient versions of
   the operators $ L^*_{\ve{\xi}}$ and $ L_{\ve{x}}$
   will be denoted by $ \bar L^*_{\ve{\xi}}$ and $\tilde L_{\ve{x}}$, respectively.
   Furthermore, we extend the resulting
   equations to $\R^3$ and denote their solutions by $\bar g$ and~$\tilde g$.
   So we get
%
   \begin{align}
     \bar L^*_{\ve\xi}\,\bar g(\ve{x};\ve\xi)
         =-\eps\,\laplace_{\ve\xi}\bar g(\ve{x};\ve\xi)
         +a(\ve{x})\,\pt_{\xi_1}\bar g(\ve{x};\ve\xi)
        &=\delta(\ve{x}-\ve\xi)\quad \mbox{for}\;\ve\xi\in\R^3,\label{eq:Green_adj_const}\\
     \widetilde L_{\ve{x}}\,\tilde g(\ve{x};\ve\xi)
          =-\eps\,\laplace_{\ve{x}}\tilde g(\ve{x};\ve\xi)
          -a(\ve\xi)\,\pt_{x_1}\tilde g(\ve{x};\ve\xi)
         &=\delta(\ve{x}-\ve\xi)\quad\mbox{for}\; \ve{x}\in\R^3.\label{eq:Green_prim_const}
   \end{align}
   As $\ve{x}$ appears in \eqref{eq:Green_adj_const} as a parameter, so
   the coefficient $a(\ve{x})$ in this equation is considered constant
   and we can solve the problem explicitly.
   Setting $q=\frac{1}{2}a(\ve{x})$ for fixed $\ve{x}\in(0,1)^3$
   and $\bar g(\ve{x};\ve{\xi}) = V(\ve{x};\ve{\xi})\,\E^{q\xi_1/\eps}$ (see, e.g., \cite{KSh87}),
   one gets
   \begin{align*}
         -\eps^2\laplace_{\ve\xi}V+q^2 V
         &=\eps \,\E^{-q\xi_1/\eps}\,\delta(\ve{x}-\ve\xi)
          =\eps \,\E^{-qx_1/\eps}\,\delta(\ve{x}-\ve\xi).
   \end{align*}
   As the fundamental solution for the operator
   $-\eps^2\laplace_{\ve\xi}+q^2$
   is $\frac{1}{4\pi\eps^2} \frac{\E^{-qr/\eps}}{r}$ 
   \cite[Chapter~VII]{TikhSamars},
   so
    \[
       V(\ve{x};\ve{\xi})
         =\eps \E^{-x_1q/\eps}\frac{1}{4\pi\eps^2}\frac{\E^{-rq/\eps}}{r}
             \quad\mbox{where}\quad
       r =\sqrt{(x_1-\xi_1)^2+(x_2-\xi_2)^2+(x_3-\xi_3)^2}.
    \]
   Finally, for the solution of \eqref{eq:Green_adj_const} we get
   \[
      \bar g(\ve{x};\ve\xi)
       =\frac{1}{4\pi\eps^2}\,\frac{\E^{q(\xi_1-x_1-r)/\eps}}{r},
              \qquad\mbox{where}\quad
       q=q(\ve{x})={\textstyle\frac{1}{2}}a(\ve{x}).
   \]
   A similar argument yields the solution of \eqref{eq:Green_prim_const}
   \[
      \tilde g(\ve{x};\ve\xi)
      =\frac{1}{4\pi\eps^2}\,\frac{\E^{q(\xi_1-x_1-r)/\eps}}{r},
              \qquad\mbox{where}\quad
       q=q(\ve\xi)={\textstyle\frac{1}{2}}a(\ve\xi).
   \]
Let $\widehat\xi_{1,[x_1]}=(\xi_1-x_1)/\eps$,
   $\widehat\xi_2=(\xi_2-x_2)/\eps$, $\widehat\xi_3=(\xi_3-x_3)/\eps$
   and $\widehat{r}_{[x_1]}=\sqrt{\widehat\xi_{1,[x_1]}^2+\widehat\xi_2^2+\widehat\xi_3^2}$.
As we shall need bounds for both $\bar g$ and $\tilde g$,
it is convenient to represent them via a more general function
   \begin{gather}\label{eq:def_g0}
     g=g(\ve{x};\ve\xi;q)
      :=\frac{1}{4\pi\eps^2}\, \frac{\E^{q(\widehat\xi_{1,[x_1]}-\widehat r_{[x_1]})}}{\widehat r_{[x_1]}}
   \end{gather}
   as
   \begin{gather}\label{bar_tilde_g_def}
      \bar g(\ve{x};\ve\xi)=g(\ve{x};\ve\xi;q)\Bigr|_{q={\textstyle\frac{1}{2}}a(\ve{x})},
      \qquad
      \tilde g(\ve{x};\ve\xi)=g(\ve{x};\ve\xi;q)\Bigr|_{q={\textstyle\frac{1}{2}}a(\ve\xi)}
   \end{gather}
   We use the subindex $[x_1]$ in $\widehat\xi_{1,[x_1]}$ and $\widehat r_{[x_1]}$
   to highlight their dependence  on $x_1$
   as in many places $x_1$ will take different values;
   but when there is no ambiguity, we shall sometimes 
   simply write
   $\widehat\xi_1$ and $\widehat r$.
%
\section{Bounds for the fundamental solution $g(\ve{x};\ve{\xi};q)$}\label{sec:bounds_green_general}
%
   Throughout this section we assume that $\Omega=(0,1)\times\R^2$,
   but all results remain valid for $\Omega=(0,1)^3$.
   Here we derive a number of useful bounds for the fundamental solution $g$
   of \eqref{eq:def_g0}
   and its derivatives that will be used in Section~\ref{sec:bounds_green_const}.
   As in $\bar g$ and $\tilde g$ we set $q=\frac12 a(\ve{x})$ and $q=\frac12 a(\ve{\xi})$,
   respectively, so  we shall also use, for $k=2,3$, the differential operators
   \begin{gather}
      D_{\xi_k}:=\pt_{\xi_k}+{\ts\frac12}\pt_{\xi_k} a(\ve{\xi})\cdot\pt_q,
         \qquad\quad
      D_{x_k}:=\pt_{x_k}+{\ts\frac12}\pt_{x_k} a(\ve{x})\cdot\pt_q.\label{D_ops}
   \end{gather}

   \begin{lem}\label{lem:g0_bounds}
      Let $\ve{x}\in[-1,1]\times\mathbb{R}^2$ and $0<\frac12\alpha\le q\le C$.
      Then for the function $g=g(\ve{x};\ve{\xi};q)$ of \eqref{eq:def_g0}
      we have the following bounds
      \begin{subequations}\label{g_bounds}
         \begin{align}
            \label{g_L1}      \norm{g(\ve{x};\cdot;q)}{1\,;\Omega}
                                 &\leq C,\\
            \label{g_xi_L1}   \norm{\pt_{\xi_1} g(\ve{x};\cdot;q)}{1\,;\Omega}
                                 &\leq C(1+|\ln\eps|),\\
            \label{g_eta_L1}  \eps^{1/2}\,\norm{\pt_{\xi_k} g(\ve{x};\cdot;q)}{1\,;\Omega}+
                                          \norm{\pt_q g(\ve{x};\cdot;q)}{1\,;\Omega}
                                 &\leq C,\quad k=2,\,3,\\
            \label{g_xi_R_L1} \norm{(\eps \widehat r_{[x_1]}\,\pt_{\xi_1} g)(\ve{x};\cdot;q)}{1\,;\Omega}
                                 &\leq C,\\
            \label{g_eta_x_L1} \eps^{1/2}\,\norm{(\eps \widehat r_{[x_1]}\,\pt^2_{\xi_1\xi_k} g)(\ve{x};\cdot;q)}{1\,;\Omega}+
                              \norm{(\eps \widehat r_{[x_1]}\,\pt^2_{\xi_1 q} g)(\ve{x};\cdot;q)}{1\,;\Omega}
                                 &\leq C,\quad k=2,\,3,
      \intertext{and for any ball $B(\ve{x}';\rho)$ of radius $\rho$ centered at any $\ve{x}'\in[0,1]\times\mathbb{R}^2$, we have}
            \label{g_eta_ball} \norm{g(\ve{x};\cdot;q)}{1,1\,;\Omega\cap B(\ve{x}';\rho)}
                                 &\leq C\eps^{-1}\rho,
      \end{align}
      while for the ball $B(\ve{x};\rho)$ of radius $\rho$ centered at $\ve{x}$, we have
      \begin{align}
         \norm{\pt^2_{\xi_1} g(\ve{x};\cdot;q)}{1\,;\Omega\setminus B(\ve{x};\rho)}
           &\leq C\eps^{-1}\ln(2+\eps/\rho),\label{g_xi2_L1}\\
         \norm{\pt^2_{\xi_k} g(\ve{x};\cdot;q)}{1\,;\Omega\setminus B(\ve{x};\rho)}
           &\leq C\eps^{-1}(\ln(2+\eps/\rho)+|\ln\eps|),\quad k=2,\,3.\label{g_eta2_L1}
      \end{align}
      \end{subequations}
      Furthermore, one has the bound
      \begin{subequations}
      \begin{gather}\label{g_x_L1}
         \norm{\pt_{x_1}   g(\ve{x};\cdot;q)}{1\,;\Omega}
            \leq C(1+|\ln\eps|),
      \end{gather}
      and  with the differential operators \eqref{D_ops}, one has, for $k=2,\,3$,
      \begin{align}
         \norm{D_{\xi_k} g(\ve{x};\cdot;q)}{1\,;\Omega}+
         \norm{D_{x_k} g(\ve{x};\cdot;q)}{1\,;\Omega}
           &\leq C\eps^{-1/2},\label{D_g}\\
         \norm{(\eps \widehat r_{[x_1]}\,D_{\xi_k}\pt_{x_1} g)(\ve{x};\cdot;q)}{1\,;\Omega}+
         \norm{(\eps \widehat r_{[x_1]}\,D_{x_k}\pt_{\xi_1} g)(\ve{x};\cdot;q)}{1\,;\Omega}
           &\leq C\eps^{-1/2}.\label{D_pt_g}
      \end{align}
      \end{subequations}
   \end{lem}

   \begin{proof}
      First, note that $\grad_{\ve{x}} g=-\grad_{\ve\xi} g$, so
      \eqref{g_x_L1} follows from \eqref{g_xi_L1},
      \eqref{D_g} follows from \eqref{D_ops},\,\eqref{g_eta_L1},
      while \eqref{D_pt_g} follows from \eqref{D_ops},\,\eqref{g_eta_x_L1}.
      Thus it suffices to establish the bounds~\eqref{g_bounds}.

       Throughout this proof,
whenever $k$ appears in any relation,
it will be understood to be valid for $k=2,3$
(as all the bounds in~\eqref{g_bounds} that involve $k$,
      are given for both $k=2,3$).

      A calculation shows that the first-order derivatives of $g=g(\ve{x};\ve{\xi};q)$
      are given by
      \begin{subequations}\label{diff_g}
      \begin{align}
      \label{g_xi1} \pt_{\xi_1} g
                     &=\frac{1}{4\pi\eps^3}\,\widehat r^{-2}\,
                           \Bigl[q(\widehat r-\widehat \xi_1)-\frac{\widehat\xi_1}{\widehat r}\,\Bigr]
                           \,\,\E^{q(\widehat\xi_1-\widehat r)},\\
      \label{g_xi2}
      \pt_{\xi_k} g
                     &=-\frac{1}{4\pi\eps^3}\,
                           \bigl(q\widehat r+1\bigr)\,\frac{\widehat\xi_k}{\widehat r^3}
                           \,\,\E^{q(\widehat\xi_1-\widehat r)},
                           \\
      \label{g_q}   \pt_q g
                     &=\frac{1}{4\pi\eps^2}\,
                           \frac{\widehat\xi_1-\widehat r}{\widehat r}
                           \,\,\E^{q(\widehat\xi_1-\widehat r)}.
      \end{align}
      \end{subequations}
      Here we used $\partial_{\xi_j}\widehat r=\eps^{-1}\widehat\xi_j/\widehat r\,$
      for $j=1,\,2,\,3$.
      In a similar manner, but also using
      $\pt_{\xi_i}(\widehat\xi_j/\widehat r)
      =-\eps^{-1}\widehat\xi_i\widehat\xi_j/\widehat r^3$ with $i\neq j$,
      one gets second-order derivatives
      \begin{subequations}\label{diff2_g}
      \begin{align}
      \label{g2_xi1_xi2} \pt^2_{\xi_1\xi_k} g
                           &=\frac{1}{4\pi\eps^4}\,
                           \frac{\widehat\xi_k}{\widehat r^3}\,
                           \Bigl[q^2(\widehat\xi_1-\widehat r)+
                                 q\frac{3\widehat\xi_1-\widehat r}{\widehat r}+
                                 3\frac{\widehat\xi_1}{\widehat r^2}\,\Bigr]
                                 \,\,\E^{q(\widehat\xi_1-\widehat r)},
                                 \\
      \label{g2_xi1_q}   \pt^2_{\xi_1 q} g
                           &=\frac{1}{4\pi\eps^3}\,
                           \widehat r^{-2}\,
                           \Bigl[-q(\widehat\xi_1-\widehat r)^2+\frac{\widehat r^2-\widehat\xi_1^2}{\widehat r}\,\Bigr]
                           \,\,\E^{q(\widehat\xi_1-\widehat r)},\\
      \label{g2_xi2_xi2} \pt^2_{\xi_k} g
                           &=\frac{1}{4\pi\eps^4}\,
                           \widehat r^{-3}\,\Bigl[q^2\widehat\xi_k^2+(q\widehat r+1)
                                                  \frac{3\widehat\xi_k^2-\widehat r^2}
                                                       {\widehat r^2}\,\Bigr]
                                                       \,\,\E^{q(\widehat\xi_1-\widehat r)}.
      \intertext{Finally, combining $\pt_{\xi_1}^2 g=-\pt_{\xi_2}^2 g-\pt_{\xi_3}^2 g+\frac{2q}\eps\,\pt_{\xi_1} g$
                 with \eqref{g_xi1} and \eqref{g2_xi2_xi2} yields}
      \label{g2_xi1_xi1} \pt^2_{\xi_1} g
                           &=\frac{1}{4\pi\eps^4}
                           \,\widehat r^{-3}\,\Bigl[ q^2\bigl(\widehat r-\widehat\xi_1\bigr)^2
                                                  -q  \bigl(\widehat r-\widehat\xi_1\bigr)
                              \Bigl(1+3\frac{\widehat\xi_1}{\widehat r}\Bigr)
                                                +\frac{3\widehat\xi_1^2-\widehat r^2}
                                                      {\widehat r^2}\,\Bigr]
                                                      \,\,\E^{q(\widehat\xi_1-\widehat r)}.
      \end{align}
      \end{subequations}
      Now we proceed to estimating the above derivatives of $g$.
      Note that $d\ve{\xi}=\eps^3 d\widehat{\ve{\xi}}$,
      where $\widehat{\ve{\xi}}\in\widehat\Omega:=\eps^{-1}(-x_1,1-x_1)\times\R^2
                                          \subset(-\infty,2/\eps)\times\R^2$.
      Consider the two sub-domains
      \[
         \widehat\Omega_1:=\bigl\{\;\widehat\xi_1<1+{\ts\frac{1}{2}}\widehat r
         \;\bigr\},
         \qquad
         \widehat\Omega_2:=\bigl\{\;
         \max\{\,1,{\ts\frac{1}{2}}\widehat r\,\}<\widehat\xi_1<2/\eps\;\bigr\}.
      \]
      As  $\widehat\Omega\subset\widehat\Omega_1\cup\widehat\Omega_2$ for any $x_1\in[-1,1]$,
      it is convenient to consider integrals over these two sub-domains separately.

     (i) Consider $\widehat{\ve{\xi}}\in\widehat\Omega_1$. Then
     $\widehat\xi_1\leq 1+\frac{1}{2}\widehat r$, so
         one gets
         \begin{align}
            \eps^3\bigl[(1+\widehat r)
                        (\eps^{-1}|g|+|\pt_{\xi_1} g|+|\pt_{\xi_k} g|
                        +|\pt_q g|+|\pt^2_{\xi_1 q} g|)
                        &+\eps\widehat r|\pt^2_{\xi_1\xi_k} g|
                        \bigr]\notag\\
               &\leq C \, \widehat r^{-2}\, (1+\widehat r+\widehat r^2+\widehat r^3)
               \,\,\E^{q(\widehat\xi_1-\widehat r)}\notag\\
               &\leq C\,\widehat r^{-2}\,\, \E^{-q\widehat r/4},\label{star0}
         \end{align}
         where we combined $\E^{q\widehat\xi_1}\le \E^{q(1+\widehat r/2)}$
         with $(1+\widehat r+\widehat r^2+\widehat r^3)\le C\E^{q\widehat r/4}$.
         This immediately yields
         \begin{multline}
            \iiint_{\widehat\Omega_1}\!
            \bigl[(1+\widehat r)
                  (\eps^{-1}|g|+|\pt_{\xi_1} g|+|\pt_{\xi_k} g|
                     +|\pt_q g|+|\pt^2_{\xi_1 q} g|)
                  +\eps\widehat r|\pt^2_{\xi_1\xi_k} g|
                  \bigr]\,\bigl(\eps^3d\widehat{\ve{\xi}}\bigr)
                  \\
            \leq C\int_0^{\infty}\E^{-q\widehat r/4} \,d\widehat r
            \leq C.\label{star}
         \end{multline}
         Similarly,
         \[
            \eps^3\bigl[|\pt^2_{\xi_1} g|+|\pt^2_{\xi_k} g|
            \bigr]
               \le C\eps^{-1}\,\widehat r^{-3}\,(1+\widehat r^2)\,\,\E^{q(\widehat\xi_1-\widehat r)}
               \le C\eps^{-1}\,\widehat r^{-2}\,(\widehat r^{-1}+\widehat r)\,\,\E^{-q\widehat r/2},
         \]
         so
         \begin{gather}\label{star2}
            \iiint_{\widehat\Omega_1\setminus B(\ve{0};\widehat\rho)}\!\bigl[|\pt^2_{\xi_1} g|+|\pt^2_{\xi_k} g|\bigr]
               \,\bigl(\eps^3d\widehat{\ve{\xi}}\bigr)
            \leq C\eps^{-1}\!\int_{\widehat\rho}^{\infty}\!\!(\widehat r^{-1}+\widehat r)\,\E^{-q\widehat r/2} \,d\widehat r
            \leq C\eps^{-1}\ln(2+\widehat\rho^{-1}).
         \end{gather}
%
         Furthermore, for an arbitrary ball $\widehat B_{\widehat\rho}$ of radius $\widehat\rho$
         in the coordinates $\ve{\widehat\xi}$, we get
         \begin{gather}\label{ball_Omega1}
            \iiint_{\widehat\Omega_1\cap \widehat B_{\widehat\rho}}\!
               \bigl[|g|+|\pt_{\xi_1} g|+|\pt_{\xi_k} g|
               \bigr]
               \,\bigl(\eps^3d\widehat{\ve{\xi}}\bigr)
            \leq C\int_0^{\widehat\rho}\!\!\E^{-q\widehat r/4} \,d\widehat r\leq C\min\{\widehat\rho,1\}.
         \end{gather}

    (ii) Next consider $\widehat{\ve{\xi}}\in\widehat\Omega_2$.
         In this sub-domain, it is convenient to rewrite the integrals in terms of
         $(\widehat \xi_1,t_2,t_3)$, where 
         \begin{gather}\label{t_def}
               t_k:=\widehat\xi_1^{-1/2}\,\widehat\xi_k,
               \quad\mbox{so}\quad
               \widehat\xi_1^{-1/2}\,d\widehat\xi_k=dt_k
               \quad\mbox{and}\quad
               \widehat r-\widehat\xi_1=\frac{\widehat\xi_2^2+\widehat\xi_3^2}{\widehat r+\widehat\xi_1}\le t_2^2+t_3^2=:t^2.
         \end{gather}
         Note that $\widehat\xi_1\le\widehat r\le 2\,\widehat\xi_1$ in $\widehat\Omega_2$
         so $\widehat{r}-\widehat\xi_1=(\widehat\xi_2^2+\widehat\xi_3^2)/(\widehat{r}+\widehat\xi_1)\ge c_0 t^2$,
         where $c_0:=\frac13
         $.
         Consequently $\E^{-q(\widehat{r}-\widehat\xi_1)}\leq \E^{-q
         c_0t^2}$ or
         \begin{gather}\label{Q_def}
         \E^{-q(\widehat{r}-\widehat\xi_1)}\leq C\, \widehat r\,
         Q,
         \qquad\mbox{where}\quad
            Q:=\widehat\xi_1^{-1}\, \E^{-q c_0t^2},
         \end{gather}
         and
         \begin{align}\notag
            \iint_{\R^2}(1+t+t^2+t^3+t^4)\,Q\,\,d\widehat\xi_2\,d\widehat\xi_3
            &= \iint_{\R^2} (1+t+t^2+t^3+t^4)\,\E^{-qc_0t^2}
               \,dt_2\,dt_3
               \\\label{Q_int}
               &\leq C.
         \end{align}
         Using \eqref{diff_g},\eqref{diff2_g} and \eqref{t_def} it is straightforward to prove the
         following bounds for $g$ and its derivatives in $\widehat\Omega_2$
         \begin{subequations}\label{bounds15}
         \begin{align}
            \label{bound0}  \eps^3| g|
                              &\leq C\,\eps\, Q,\\
            \label{bound1}  \eps^3|\pt_{\xi_k} g|
                              &\leq C\,\widehat\xi_1^{-1/2}\, t\, Q,
                              \\
            \label{bound2}  \eps^3|\pt^2_{\xi_k} g|
                              &\leq C\,\eps^{-1}\,\widehat\xi_1^{-1}\, [1+t^2]\, Q,
%
\intertext{and also}
            \label{bound3}  \eps^3(\eps \widehat r|\pt_{\xi_1} g|+|\pt_q g|)
                              &\leq C\,\eps\, [1+t^2]\, Q,\\
            \label{bound4}  \eps^3|\pt_{\xi_1} g|
                              &\leq C\,\widehat\xi_1^{-1}\, [1+t^2]\, Q,\\
            \label{bound5}  \eps^3(\eps\widehat r|\pt_{\xi_1\xi_k} g|)
                              &\leq C\,\widehat\xi_1^{-1/2}\, t\,[1+t^2]\, Q,
                              \\
            \label{bound6}  \eps^3(\eps \widehat r|\pt^2_{\xi_1 q} g|)
                              &\leq C\,\eps\, (t^2+t^4)\, Q,\\
            \label{bound7}  \eps^3|\pt^2_{\xi_1} g|
                              &\leq C\,\eps^{-1}\, \widehat\xi_1^{-2}\,(1+t^2+t^4)\, Q.
         \end{align}
         \end{subequations}

         Combining the obtained estimates \eqref{bounds15} with \eqref{Q_int} yields
         \begin{align}
         \hspace{-0.2cm}
            \iiint_{\widehat\Omega_2}\!
             \bigl[|g|+\eps^{1/2}|\pt_{\xi_k} g|
             +\eps\widehat r|\pt_{\xi_1} g|+|\pt_q g|
             &+\eps^{1/2}\eps\widehat r|\pt^2_{\xi_1\xi_k} g|
             +\eps\widehat r|\pt^2_{\xi_1 q} g|
             +|\pt_{\xi_1}^2 g|\bigr]
             \,\bigl(\eps^3d\widehat{\ve{\xi}}\bigr)\notag\\
             &\leq C\int_1^{2/\eps}\!\![\eps+\eps^{1/2}\widehat\xi_1^{-1/2}] \,
             d\widehat\xi_1
             \leq C.\label{int_Omega2_main}
         \end{align}
         Similarly, combining  \eqref{bound2} and \eqref{bound4} with \eqref{Q_int} yields
         \begin{gather}\label{int_Omega2_aux}
            \iiint_{\widehat\Omega_2}\!\bigl[|\pt_{\xi_1} g|
            +\eps|\pt^2_{\xi_k} g|\bigr]
               \,\bigl(\eps^3d\widehat{\ve{\xi}}\bigr)
            \leq C\int_1^{2/\eps}\!\!\widehat\xi_1^{-1} \,d\widehat\xi_1
            \leq C(1+|\ln\eps|).
         \end{gather}
         Furthermore, by \eqref{bound1}, and \eqref{bound4}
         for an arbitrary ball $\widehat B_{\widehat\rho}$ of radius $\widehat\rho$
         in the coordinates $\widehat{\ve{\xi}}$, we get
         \begin{gather}\label{ball_Omega2}
            \iiint_{\widehat\Omega_2\cap \widehat B_{\widehat\rho}}\!
              (|g|+|\pt_{\xi_1} g|
              +|\pt_{\xi_k} g|)
                    \,\bigl(\eps^3d\widehat{\ve{\xi}}\bigr)
            \leq C\int_1^{1+\widehat\rho}\!\!
            \bigl[\eps+\widehat\xi_1^{-1}+\widehat\xi_1^{-1/2}\bigr]
            \,d\widehat\xi_1
            \leq C\widehat\rho.
         \end{gather}

      To complete the proof, we now recall that $\widehat\Omega\subset\widehat\Omega_1\cup\widehat\Omega_2$ and
      combine estimates \eqref{star} and \eqref{star2} (that involve integration over $\widehat\Omega_1$)
      with \eqref{int_Omega2_main} and \eqref{int_Omega2_aux},
      which yields the desired bounds \eqref{g_L1}-\eqref{g_eta_x_L1}
      and \eqref{g_xi2_L1}, \eqref{g_eta2_L1}.
      To get the latter two bound we also used the observation that the ball $B(\ve{x};\rho)$ of radius $\rho$
      in the coordinates $\ve{\xi}$ becomes the ball $B(\ve{0};\widehat\rho)$ of radius $\widehat\rho=\eps^{-1}\rho$
      in the coordinates $\widehat{\ve{\xi}}$.
      The remaining assertion \eqref{g_eta_ball}
      is obtained by combining \eqref{ball_Omega1} with \eqref{ball_Omega2} and
      noting that an arbitrary ball $B(\ve{x}';\rho)$ of radius $\rho$
      in the coordinates $\ve{\xi}$
      becomes a ball $\widehat B_{\widehat\rho}$ of radius $\widehat\rho=\eps^{-1}\rho$
      in the coordinates $\widehat{\ve{\xi}}$.
   \end{proof}

   Our next result shows that for $x_1\ge 1$, one gets stronger bounds for $g$ and its derivatives.
   These bounds involve the weight function
   \begin{gather}\label{lambda_def}
      \lambda:=\E^{2q(x_1-1)/\eps}.
   \end{gather}
   and show that, although $\lambda$ is exponentially large in $\eps$,
   this is compensated by the smallness of $g$ and its derivatives.
   \begin{lem}\label{lem:g_lmbd_bounds}
      Let $\ve{x}\in[1,3]\times\R^2$ and $0<\frac12\alpha\le q\le C$.
      Then for the function $g=g(\ve{x};\ve{\xi};q)$ of \eqref{eq:def_g0}
      and the weight $\lambda$ of \eqref{lambda_def},
      one has the following bounds
      \begin{subequations}\label{g_bounds_lmbd_}
      \begin{align}
         \norm{([1+\eps\widehat r_{[x_1]}]\,\lambda g)(\ve{x};\cdot;q)}{1\,;\Omega}
            &\leq C\eps,\label{g_L1_mu}\\
         \norm{(\lambda\,\pt_{\xi_1} g)(\ve{x};\cdot;q)}{1\,;\Omega}
        +\norm{(\lambda\,\pt_q g)(\ve{x};\cdot;q)}{1\,;\Omega}
            &\leq C,\label{g_L1_mu2}\\
         \norm{([1+\eps^{1/2}\widehat r_{[x_1]}]\,\lambda\,\pt_{\xi_k} g)(\ve{x};\cdot;q)}{1\,;\Omega}
        +\eps^{1/2}\norm{(\eps\widehat r_{[x_1]}\,\lambda\,\pt_{\xi_1\xi_k}^2 g)(\ve{x};\cdot;q)}{1\,;\Omega}
            &\leq C,\quad k=2,3,
            \label{g_L1_mu3}\\
         \norm{\widehat r_{[x_1]}\,\pt_q(\lambda\, g)(\ve{x};\cdot;q)}{1\,;\Omega}
        +\norm{\eps \widehat r_{[x_1]}\,\pt_q(\lambda\,\pt_{\xi_1} g)(\ve{x};\cdot;q)}{1\,;\Omega}
            &\le C,\label{g_L1_mu4}
      \intertext{and for any ball $B(\ve{x}';\rho)$ of radius $\rho$ centered at any
                 $\ve{x}'\in[0,1]\times\R^2$, one has}
        \norm{(\lambda\, g)(\ve{x};\cdot;q)}{1,1\,;\Omega\cap B(\ve{x}';\rho)}
            &\leq C\eps^{-1}\rho,\label{g_L1_mu5}
      \end{align}
      while for the ball $B(\ve{x};\rho)$ of radius $\rho$ centered at $\ve{x}$ and $k=1,\,2,\,3$, one has
      \begin{gather}\label{g_L1_mu6}
          \norm{(\lambda\,\pt^2_{\xi_k} g)(\ve{x};\cdot;q)}{1\,;\Omega\setminus B(\ve{x},\rho)}
               \leq C\eps^{-1}\ln(2+\eps/\rho).
      \end{gather}
      \end{subequations}
      Furthermore, with the differential operators \eqref{D_ops} and $k=2,\,3$, we have
      \begin{subequations}\label{g_lmb_bounds}
         \begin{align}
             \norm{\pt_{x_1}(\lambda  g)(\ve{x};\cdot;q)}{1\,;\Omega}
            +\norm{D_{x_k}(\lambda  g)(\ve{x};\cdot;q)}{1\,;\Omega}
            +\norm{D_{\xi_k}(\lambda  g)(\ve{x};\cdot;q)}{1\,;\Omega}
            &\leq C,\label{D_g_lmbd}\\
             \norm{\eps \widehat r_{[x_1]}\,D_{x_k}(\lambda\, \pt_{\xi_1} g)(\ve{x};\cdot;q)}{1\,;\Omega}
            +\norm{\eps \widehat r_{[x_1]}\,D_{\xi_k}\pt_{x_1}(\lambda g)(\ve{x};\cdot;q)}{1\,;\Omega}
            &\leq C\eps^{-1/2}.\label{D_pt_lmbd}
         \end{align}
      \end{subequations}
   \end{lem}
   \begin{proof}
      Throughout this proof,
whenever $k$ appears in any relation, it will be understood to be
valid for $k=2,3$ (as all the bounds in~\eqref{g_bounds_lmbd_},
\eqref{g_lmb_bounds} that involve $k$,
      are given for both $k=2,3$).

      We  shall use the notation $A=A(x_1):=(x_1-1)/\eps\geq 0$.
      Then \eqref{lambda_def} becomes $\lambda=\E^{2qA}$.
      We partially imitate the proof of Lemma~\ref{lem:g0_bounds}.
      Again $d\ve{\xi}=\eps^3\,d\widehat{\ve{\xi}}$,
      but now  $\widehat{\ve{\xi}}\in\widehat\Omega=\eps^{-1}(-x_1,1-x_1)\times\R^2
                \subset(-3/\eps,-A)\times\R^2$.
      So $\widehat\xi_1<-A\le 0$ immediately yields
      \begin{gather}\label{lambda_exp}
         \lambda\, \E^{q \widehat\xi_1}=\E^{2q(A-|\widehat\xi_1|)}\,\E^{q|\widehat\xi_1|}\le \E^{q|\widehat\xi_1|}.
      \end{gather}

      Consider the sub-domains
      \begin{align*}
         \widehat\Omega'_1&:=\bigl\{\;|\widehat\xi_1|<1+{\ts\frac12\widehat r},\;\;\widehat\xi_1<-A\;\bigr\},\\
         \widehat\Omega'_2&:=\bigl\{\;|\widehat\xi_1|>\max\{\ts1,{\ts\frac12\widehat r}\},\;\;
                                    -3/\eps<\widehat\xi_1<-A\;\bigr\}.
      \end{align*}
      As $\widehat\Omega\subset\widehat\Omega_1'\cup\widehat\Omega_2'$
      for any $x_1\in[1,3]$, we estimate integrals over these two domains
      separately.

  (i) Let $\widehat{\ve{\xi}}\in\widehat\Omega'_1$. Then $|\widehat\xi_1|\leq 1+\frac12\widehat r$
      so, by \eqref{lambda_exp}, one has $\lambda\, \E^{q \widehat\xi_1}\le \E^{q(1+\widehat r/2)}$.
      The first inequality in \eqref{star0} remains valid, but now we combine it with
      \begin{gather}\label{eq:lambda_eqxi}
         \lambda \,\E^{q(\widehat\xi_1-\widehat r)}\, (1+\widehat r+\widehat r^2+\widehat r^3)\,
         \le C\, \E^{-q\widehat r/4}
      \end{gather}
      (which is obtained similarly to the final line in \eqref{star0}).
      This leads to a version of \eqref{star} that involves the weight $\lambda$:
      \begin{multline}\label{star_lambda}
         \iiint_{\widehat\Omega_1'}\lambda\,
         \bigl[(1+\widehat r)
               (\eps^{-1}|g|+|\pt_{\xi_1} g|+|\pt_{\xi_k} g|
                  +\eps^{-1}|\pt_q g|+|\pt^2_{\xi_1 q} g|)
               +\eps\widehat r|\pt^2_{\xi_1\xi_k} g|\bigr]
         \,\bigl(\eps^3d\widehat{\ve{\xi}}\bigr)
         \leq C.
      \end{multline}
      In a similar manner, we obtain versions of estimates \eqref{star2}
      and \eqref{ball_Omega1}, that also involve the weight $\lambda$:

      \begin{gather}\label{star2_lambda}
       \iiint_{\widehat\Omega_1'\setminus B(\ve{0};\widehat\rho)}\!
        \lambda\,|\pt^2_{\xi_k} g|
         \,\bigl(\eps^3d\widehat{\ve{\xi}}\bigr)
       \leq C\eps^{-1}\ln( 2+\widehat\rho^{-1}),
      \end{gather}
      \begin{gather}\label{ball_Omega1_lambda}
        \iiint_{\widehat\Omega_1'\cap \widehat B_{\widehat\rho}}\!
          \lambda\bigl[|g|+|\pt_{\xi_1} g|
          +|\pt_{\xi_k} g|\bigr]
        \,\bigl(\eps^3d\widehat{\ve{\xi}}\bigr)
        \leq C\min\{\widehat\rho,1\},
      \end{gather}
      where $\widehat B_{\widehat\rho}$ is  an arbitrary ball of radius $\widehat\rho$
      in the coordinates $\widehat{\ve{\xi}}$.
      Furthermore, \eqref{star_lambda} combined with
      $|\pt_q(\lambda\, g)|\le \lambda(2A|g|+|\pt_q g|)$
      and $|\pt_q(\lambda\,\pt_{\xi_1} g)|\leq \lambda(2A|\pt_{\xi_1} g|+|\pt^2_{\xi_1 q} g|)$
      and then with $A\le 2/\eps$ yields
      \begin{gather}\label{star_lambda2}
       \iiint_{\widehat\Omega_1'}\, \widehat r\,\bigl[
              |\pt_q(\lambda\, g)|+\eps|\pt_q(\lambda\,\pt_{\xi_1} g)|\bigr]
        \,\bigl(\eps^3d\widehat{\ve{\xi}}\bigr)\leq C.
      \end{gather}

 (ii) Now consider $\widehat{\ve{\xi}}\in\widehat\Omega'_2$.
      In this sub-domain (similarly to $\widehat\Omega_2$ in the proof of Lemma~\ref{lem:g0_bounds})
      one has $|\widehat\xi_1|\le \widehat r\le2|\widehat\xi_1|$ and
      $c_0 t^2\le\widehat r-|\widehat\xi_1|\le t^2$, where
      $t_k:=|\widehat\xi_1|^{-1/2}\,\widehat\xi_k$ for $k=2,3$,
      and
      $t^2:=t_2^2+t_3^2$, (compare with \eqref{t_def}).
      We also introduce a new barrier $Q$
      \begin{gather}\label{Q_def2}
         Q:=\lambda^{-1}\,\E^{2q(A-|\widehat\xi_1|)}\,\bigl\{|\widehat\xi_1|^{-1}\, \E^{-q c_0t^2}\bigr\}
         \qquad\Rightarrow\quad
          \E^{-q(\widehat{r}-\widehat\xi_1)}\leq C\, \widehat r\,
         Q,
      \end{gather}
      (compare with \eqref{Q_def}; to get the bound for $\E^{-q(\widehat{r}-\widehat\xi_1)}$ we used
      \eqref{lambda_exp}).

      With the new definition \eqref{Q_def2} of $Q$, the bounds \eqref{bound0}--\eqref{bound2}
      remain valid in $\widehat\Omega_2'$ only with $\widehat\xi_1$ replaced by $|\widehat\xi_1|$.
      Note that the bounds \eqref{bound3}--\eqref{bound6} are not valid in $\widehat\Omega_2'$,
      (as they were obtained using $\widehat r-\widehat\xi_1\leq t^2$,
       which is not the case for $\widehat\xi_1<0$).
      Instead, using $\widehat r\geq |\widehat\xi_1|\geq 1$ and
      $\widehat r\leq2|\widehat\xi_1|$,
      we prove, directly from \eqref{diff_g},\eqref{diff2_g},
      the following bounds in $\widehat\Omega'_2$:
      \begin{subequations}\label{bounds_lambda}
         \begin{align}
            \eps^3|\pt_{\xi_1} g|
               &\leq C\, Q, \label{bound2b_lmb}\\
            \eps^3|\pt_q g|
               &\leq C\,\eps |\widehat\xi_1|\, Q,\label{bound2_lmb}\\
            \eps^3(\eps\widehat r|\pt_{\xi_1\xi_k} g|)
               &\leq  C\,|\widehat\xi_1|^{1/2}\, t\, Q,
               \label{bound5_lmb}\\
            \eps^3(|\pt_q(\lambda\, g)|+\eps|\pt_q(\lambda\,\pt_{\xi_1} g)|)
               &\leq C\,\eps \lambda\, [(|\widehat\xi_1|-A)+t^2+1]\, Q.\label{bound2a_lmb}
         \end{align}
      \end{subequations}
      In particular, to establish \eqref{bound2a_lmb}, we combined
      $\pt_q(\lambda\, g)=\lambda[2A\,g+\pt_{q}g]$ and
      $\pt_q(\lambda\,\pt_{\xi_1} g)=\lambda[2A\,\pt_{\xi_1} g+\pt^2_{\xi_1 q}g]$
      with the observations that
      \[
         (\widehat{r}+|\widehat\xi_1|)-2A
            = 2(|\widehat\xi_1|-A)+(\widehat r-|\widehat\xi_1|)
            \leq 2(|\widehat\xi_1|-A)+t^2
      \]
      and $\widehat r^{-1}A\leq C$.

      Next, note that \eqref{Q_int} is valid with $Q$
      replaced by
 the multiplier $\bigl\{|\widehat\xi_1|^{-1}\, \E^{-q c_0t^2}\bigr\}$ from
      the current definition \eqref{Q_def2}
      of $Q$.
      Combining this observation with the bounds \eqref{bound0}--\eqref{bound2} and
      \eqref{bound2b_lmb}--\eqref{bound5_lmb}, and also with
      $\widehat r\le 2|\widehat\xi_1|$, yields
      \begin{align}
         \iiint_{\widehat\Omega_2'}\!\! \lambda\,
         \bigl[&(\eps^{-1}+\widehat r)|g|+|\pt_{\xi_1} g|
                  +(1+\eps^{1/2}\widehat r)|\pt_{\xi_k} g|
                  +|\pt_q g|
               +\eps^{1/2}(\eps\widehat r|\pt^2_{\xi_1\xi_k} g|)
                  +\eps|\pt^2_{\xi_k} g|\bigr]
         \,\bigl(\eps^3d\widehat{\ve{\xi}}\bigr)\notag\\
         &\leq C\int_{-3/\eps}^{-\max\{A,1\}}\!
            \bigl[1+\eps|\widehat\xi_1|+|\widehat\xi_1|^{-1/2}
                  +(\eps|\widehat\xi_1|)^{1/2}+|\widehat\xi_1|^{-1}\bigr]\,
            \E^{2q(A-|\widehat\xi_1|)} \,d\widehat\xi_1
          \leq C.\label{int_Omega2_main_lmbd}
      \end{align}
      Similarly, from \eqref{bound2a_lmb} combined with
      $\widehat r\leq 2|\widehat\xi_1|\le 6\eps^{-1}$, one gets
      \begin{multline}\label{Omega2_prime}
         \iiint_{\widehat\Omega_2'}\!
           \widehat r\bigl[|\pt_q(\lambda\, g)|+\eps|\pt_q(\lambda\,\pt_{\xi_1} g)|\bigr]
            \,\bigl(\eps^3d\widehat{\ve{\xi}}\bigr)\\
         \leq C\int_{-3/\eps}^{-\max\{A,1\}}\!\!
               \bigl[(|\widehat\xi_1|-A)+1\bigr]\,\E^{2q(A-|\widehat\xi_1|)} \,d\widehat\xi_1
         \leq C.
      \end{multline}
      Furthermore, by \eqref{bound1}, and \eqref{bound2b_lmb},
      for an arbitrary ball $\widehat B_{\widehat\rho}$ of radius $\widehat\rho$
      in the coordinates $\widehat{\ve{\xi}}$, we get
      \begin{align}\notag
         \iiint_{\widehat\Omega'_2\cap \widehat B_{\widehat\rho}}\!
            \lambda\,[|g|+|\pt_{\xi_1} g|
                      +|\pt_{\xi_k} g|]\,\bigl(\eps^3d\widehat{\ve{\xi}}\bigr)
         &\leq C\int^{-\max\{A,1\}}_{-\max\{A,1\}-\widehat\rho}
            \bigl[1+|\widehat\xi_1|^{-1/2}\bigr]\,\E^{2q(A-|\widehat\xi_1|)}d\widehat\xi_1
\\\label{ball_Omega2_lambda}
         &\leq C\widehat\rho.
      \end{align}

      To complete the proof of \eqref{g_bounds_lmbd_}, we now
      recall that $\widehat\Omega\subset\widehat\Omega_1'\cup\widehat\Omega_2'$ and
      combine  estimates \eqref{star_lambda}, \eqref{star2_lambda}, \eqref{star_lambda2}
      (that involve integration over $\widehat\Omega_1'$)
      with \eqref{int_Omega2_main_lmbd}, \eqref{Omega2_prime},
      which yields the desired bounds \eqref{g_L1_mu}--\eqref{g_L1_mu4}
      and the bounds for $\pt_{\xi_2}^2g$ and $\pt_{\xi_3}^2 g$ in \eqref{g_L1_mu6}.
      To get the latter two bounds we also used the observation
      that the ball $B(\ve{x};\rho)$ of radius $\rho$
      in the coordinates $\ve{\xi}$
      becomes the ball $B(\ve{0};\widehat\rho)$ of radius $\widehat\rho=\eps^{-1}\rho$
      in the coordinates $\widehat{\ve{\xi}}$.
      The bound for $\pt^2_{\xi_1} g$ in \eqref{g_L1_mu6}
      follows as $\pt^2_{\xi_1} g=-\pt^2_{\xi_2} g-\pt^2_{\xi_3} g+\frac{2q}\eps\,\pt_{\xi_1} g$
      for $\ve{\xi}\neq\ve{x}$.
      The remaining assertion \eqref{g_L1_mu5} is obtained by combining
      \eqref{ball_Omega1_lambda} with \eqref{ball_Omega2_lambda} and
      noting that an arbitrary ball $B(\ve{x}';\rho)$ of radius $\rho$
      in the coordinates $\ve{\xi}$
      becomes a ball $\widehat B_{\widehat\rho}$ of radius $\widehat\rho=\eps^{-1}\rho$
      in the coordinates $\widehat{\ve{\xi}}$.
      Thus we have established all the bounds \eqref{g_bounds_lmbd_}.

      We now proceed to the proof of the bounds \eqref{g_lmb_bounds}.
      Note that $\grad_{\ve{x}} g=-\grad_{\ve\xi} g$.
      Combining these with \eqref{g_L1_mu2} and
      the bounds for $\norm{\lambda\,\pt_{\xi_2} g}{1\,;\Omega}$
      and $\norm{\lambda\,\pt_{\xi_3} g}{1\,;\Omega}$ in \eqref{g_L1_mu3},
      yields
      \[
          \norm{\lambda\, \pt_{x_1} g}{1\,;\Omega}
         +\norm{\lambda\, D_{x_k} g}{1\,;\Omega}
         +\norm{\lambda\, D_{\xi_k} g}{1\,;\Omega}
         \leq C.
      \]
      Now, combining
      $\pt_{x_1}\lambda=2q\eps^{-1}\lambda$ and
      $\pt_q\lambda=2A\lambda\le 4\eps^{-1}\lambda$
      with \eqref{g_L1_mu}, yields
      \[
          \norm{g\,\pt_{x_1}\lambda}{1\,;\Omega}
         +\norm{g\,D_{x_k}\lambda}{1\,;\Omega}
         +\norm{g\,D_{\xi_k}\lambda}{1\,;\Omega}
         \leq C.
      \]
      Consequently, we get \eqref{D_g_lmbd}.

      To estimate $\eps \widehat r\,D_{x_k}(\lambda\, \pt_{\xi_1} g)$,
      note that it involves
      $\eps \widehat r\,\pt_{x_k}(\lambda\, \pt_{\xi_1} g)
      =-\eps \widehat r\,\lambda\, \pt^2_{\xi_1\xi_k} g$
      for which we have a bound in \eqref{g_L1_mu3},
      and also $\eps \widehat r\,\pt_q(\lambda\, \pt_{\xi_1} g)$,
      for which we have a bound in \eqref{g_L1_mu4}. The desired bounds
      for $\eps \widehat r\,D_{x_k}(\lambda\, \pt_{\xi_1} g)$ 
      in \eqref{D_pt_lmbd} follow.

      For $\eps \widehat r\,D_{\xi_k}\pt_{x_1}(\lambda g)$
      in \eqref{D_pt_lmbd}, a calculation yields
      $\eps \widehat r\,D_{\xi_k}\pt_{x_1}(\lambda g)
      =\eps \widehat r\,D_{\xi_k}(\lambda\, \pt_{x_1} g)+2\widehat r\,D_{\xi_k}
      ( q \lambda\,g)$.
      The first term is estimated similarly to
      $\eps \widehat r\,D_{x_k}(\lambda\, \pt_{\xi_1} g)$ in \eqref{D_pt_lmbd}.
      The remaining term $\widehat r\,D_{\xi_k}( q \lambda\,g)$ involves
      $\widehat r\,\pt_{\xi_k}( q \lambda\,g)=q\,\widehat r\,\lambda\,\pt_{\xi_k} g$,
      for which we have a bound in \eqref{g_L1_mu3}, and also
      $\widehat r\,\pt_q( q \lambda\,g)=q \,\widehat r\,\pt_q( \lambda\,g)+\widehat r\,\lambda\,g$,
      for which we have bounds in \eqref{g_L1_mu4} and \eqref{g_L1_mu}.
      Consequently \eqref{D_pt_lmbd} is proved.
   \end{proof}

   \begin{lem}\label{lem_1_3}
      Under the conditions of Lemma~\ref{lem:g_lmbd_bounds},
      for some positive constant $c_1$ one has
      \begin{gather}\label{lambda_g_1_3}
          \|\lambda g(\ve{x};\cdot)\|_{2,1\,;[0\frac13]\times\R^2}
         +\|D_{x_k}(\lambda g)(\ve{x};\cdot)\|_{1,1\,;[0\frac13]\times\R^2}
         \le  C\E^{-c_1\alpha/\eps},\quad k=2,\,3.
      \end{gather}
   \end{lem}
   \begin{proof}
      We imitate the proof of Lemma~\ref{lem:g_lmbd_bounds}, only now
      $\xi_1<\frac13$ or $\widehat\xi_1< (\frac13-x_1)/\eps\le -\frac23/\eps$.
      Thus instead of the sub-domains $\widehat\Omega_1'$ and $\widehat\Omega_2'$
      we now consider $\widehat\Omega_1''$ and $\widehat\Omega_2''$
      defined by
      $\widehat\Omega_k'':=\widehat\Omega_k'\cap\{\widehat\xi_1< -(x_1-\frac13)/\eps\}$.
      Thus in $\widehat\Omega_1''$ \eqref{eq:lambda_eqxi} remains valid with $q\ge\frac12\alpha$,
      but now $\widehat r>\frac23/\eps$.
      Therefore, when we integrate over $\widehat\Omega_1''$
      (instead of $\widehat\Omega_1'$),
      the integrals of type \eqref{star_lambda}, \eqref{star2_lambda}
      become bounded by $C\E^{-c_1\alpha/\eps}$ for any fixed $ c_1<\frac1{8}$.
      Next, when considering integrals over $\widehat\Omega_2''$
      (instead of $\widehat\Omega_2'$),
      note that $A-|\widehat\xi_1|\le-\frac23/\eps$ so
      the quantity $\E^{2q(A-|\widehat\xi_1|)}$ in the definition \eqref{Q_def2} of $Q$
      is now bounded by $\E^{-\frac23\alpha/\eps}$.
      Consequently, the integrals of type \eqref{int_Omega2_main_lmbd}
      over $\widehat\Omega_2''$ also become bounded by
      $C\E^{-c_1\alpha/\eps}$.
   \end{proof}
   \begin{rem}\label{rem_q_var}
      The estimates of Lemmas~\ref{lem:g0_bounds} and~\ref{lem:g_lmbd_bounds} remain valid
      if we set $q:=\frac12a(\ve{x})$ or $q:=\frac12 a(\ve{\xi})$
      in $g$, $\lambda$, and their derivatives (after the differentiation is performed).
   \end{rem}

\section{Approximations $\bar G$ and $\tilde G$ for the Green's function~$G$}\label{sec:bounds_green_const}
%
   We shall use two related cut-off functions $\omega_0$ and $\omega_1$
   defined by
   \begin{gather}\label{cut_off}
      \omega_0(t) \in C^2(0,1),\quad\omega_0(t)=1\;\mbox{for }t\le{\ts\frac23},
      \quad\omega_0(t)=0\;\mbox{for }t\ge{\ts\frac56};
      \quad
      \omega_1(t):=\omega_0(1-t),
   \end{gather}
   so $\omega_m(m)=1$, $\omega_m(1-m)=0$
   and $\omega_m'(t)\bigr|_{t=0,1}=\omega_m''(t)\bigr|_{t=0,1}=0$ for $m=0,1$.

   Our purpose in this section is to introduce and estimate
   frozen-coefficient approximations
   $\bar G$ and $\tilde G$ of $G$. We consider
   the domain $\Omega=(0,1)\times\R^2$
   in the first part of this section, and the domain $\Omega=(0,1)^3$ in the second part.
   Note that although $\bar G$ and $\tilde G$ will be constructed as solution approximations
   for the frozen-coefficient equations,
   we shall see in Section~\ref{sec:main_proof} that they, in fact,
   provide approximations to the Green's function $G$ for our original
   variable-coefficient problem.
%
 \subsection{Approximations
$\bar G$ and $\tilde G$
  in the domain $\Omega=(0,1)\times\R^2$}\label{ssec:bounds_green_const_1}
To construct approximations
$\bar G$ and $\tilde G$, we employ the method of images
with an inclusion of the
cut-off functions of \eqref{cut_off}.
So, using the fundamental solution $g$ of \eqref{eq:def_g0}, we define
   \begin{gather}\label{bar_tilde_G}
      \bar G(\ve{x};\ve\xi):=\bar{\mathcal G}
      \bigr|_{q=\frac12a(\ve{x})},
      \qquad
      \tilde G(\ve{x};\ve\xi):=\tilde {\mathcal G}
      \bigr|_{q=\frac12a(\ve\xi)},
   \end{gather}
   \vspace{-0.8cm}
   \begin{subequations}\label{bar_tilde_G20}
   \begin{align}
      \bar {\mathcal G}(\ve{x};\ve\xi;q)
        &:=\frac{\E^{q\widehat\xi_{1,[x_1]}}}{4\pi\eps^2}
           \left\{
              \left[ \frac{\E^{-q\widehat r_{[  x_1]}}}{\widehat r_{[  x_1]}}
                    -\frac{\E^{-q\widehat r_{[ -x_1]}}}{\widehat r_{[ -x_1]}}\right]
             -\left[ \frac{\E^{-q\widehat r_{[2-x_1]}}}{\widehat r_{[2-x_1]}}
                    -\frac{\E^{-q\widehat r_{[2+x_1]}}}{\widehat r_{[2+x_1]}}\right]\omega_1(\xi_1)
           \right\},\label{bar_G_def}\\
      \tilde {\mathcal G}(\ve{x};\ve\xi;q)
        &:=\frac{\E^{q\widehat\xi_{1,[x_1]}}}{ 4\pi\eps^2}
           \left\{
              \left[ \frac{\E^{-q\widehat r_{[  x_1]}}}{\widehat r_{[  x_1]}}
                    -\frac{\E^{-q\widehat r_{[2-x_1]}}}{\widehat r_{[2-x_1]}}\right]
             -\left[ \frac{\E^{-q\widehat r_{[ -x_1]}}}{\widehat r_{[ -x_1]}}
                    -\frac{\E^{-q\widehat r_{[2+x_1]}}}{\widehat r_{[2+x_1]}}\right]\omega_0(x_1)
           \right\}\label{tilde_G_def}.
   \end{align}
   \end{subequations}
   Note that $\bar G\bigr|_{\xi_1=0,1}=0$ and
   $\tilde G\bigr|_{x_1=0,1}=0$
   (the former observation
   follows from $\widehat r_{[x_1]}=\widehat r_{[-x_1]}$ at $\xi_1=0$,
   and $\widehat r_{[x_1]}=\widehat r_{[2-x_1]}$ and $\widehat r_{[-x_1]}=\widehat r_{[2+x_1]}$ at $\xi_1=1$).
   We shall see shortly (see Lemma~\ref{lem_tilde_bar_w}) that
   $\bar L^*_{\ve\xi}\bar G\approx L^*_{\ve\xi}G$
   and $\tilde L_{\ve{x}}\tilde G\approx L_{\ve{x}} G$;
   in this sense $\bar G$ and $\tilde G$ give approximations for $G$.

   Rewrite the definitions of $\bar{\mathcal G}$ and $\tilde{\mathcal G}$
   using the notation
   \begin{subequations}
   \begin{align}
      g_{[d]}
        &:= g(d,x_2,x_3;\ve\xi;q)
          ={\frac{1}{ 4\pi\eps^2}}\, \frac{\E^{q(\widehat\xi_{1,[d]}-\widehat r_{[d]})}}{\widehat r_{[d]}},\label{g_x_brackets}\\
      \lambda^{\pm}&:=\E^{2q(1\pm x_1)/{\eps}},
      \qquad\quad  p:=\E^{-2qx_1/{\eps}},\label{lmb_p_def}
   \end{align}
   \end{subequations}
   and the observation that
   \begin{gather}\label{q_d}
      {\frac{1}{ 4\pi\eps^2}}\,\frac{\E^{q(\widehat\xi_{1,[x_1]}-\widehat r_{[d]})}}{\widehat r_{[d]}}
      =\E^{q(d-x_1)/\eps}\,g_{[d]}\qquad\mbox{for}\;\; d=\pm x_1, 2\pm x_1.
   \end{gather}
   They yield
   \begin{subequations}\label{bar_G_tilde_g}
   \begin{align}
      \bar{\mathcal G}(\ve{x};\ve\xi;q)
       &= \left[g_{[x_1]}-p\, g_{[-x_1]}\right]
          -\left[\lambda^- g_{[2-x_1]}-p\,\lambda^{\!+} g_{[2+x_1]}\right]\omega_1(\xi_1),\label{bar_G_g}\\
      \tilde {\mathcal G}(\ve{x};\ve\xi;q)
       &= \left[g_{[x_1]}-\lambda^- g_{[2-x_1]}\right]
         -\left[p\, g_{[-x_1]}-p\,\lambda^{\!+} g_{[2+x_1]}\right]\omega_0(x_1).\label{tilde_G_g}
   \end{align}
   \end{subequations}
   Note that $\lambda^\pm$ is obtained by replacing $x_1$ by $2\pm x_1$
   in the definition \eqref{lambda_def} of $\lambda$.

   In the next lemma, we estimate the functions
   \begin{gather}\label{tilde_w_def}
      \bar \phi(\ve{x};\ve\xi):=\bar L^*_{\ve\xi}\bar G-L^*_{\ve\xi}G,
      \qquad
      \tilde \phi(\ve{x};\ve\xi):=\tilde L_{\ve{x}}\tilde G-L_{\ve{x}}G.
   \end{gather}

   \begin{lem}\label{lem_tilde_bar_w}
      Let $\ve{x}\in\Omega=(0,1)\times\R^2$. Then
      for the functions $\bar\phi$ and $\tilde \phi$ of \eqref{tilde_w_def}, one has
      \begin{gather}\label{tilde_bar_w}
          \|\bar \phi(\ve{x};\cdot)\|_{1,1\,;\Omega}
         +\|\pt_{x_2}\bar \phi(\ve{x};\cdot)\|_{1\,;\Omega}
         +\|\pt_{x_3}\bar \phi(\ve{x};\cdot)\|_{1\,;\Omega}
         +\|\tilde \phi(\ve{x};\cdot)\|_{1,1\,;\Omega}
         \leq C\E^{-c_1\alpha/\eps}
         \leq C.
      \end{gather}
      One also has
      \begin{gather}\label{eq:bar_phi_boundary}
         \bar\phi(\ve{x};\ve{\xi})|_{\ve{\xi}\in\pt\Omega}=0.
      \end{gather}
   \end{lem}

   \begin{proof}
     (i) First we prove the desired assertions for $\bar\phi$. By \eqref{bar_tilde_G},
         throughout this part of the proof we set $q=\frac12a(\ve{x})\ge\frac12\alpha$.
         Recall that $\bar g$ solves the differential equation \eqref{eq:Green_adj_const}
         with the operator $\bar L^*_{\ve\xi}$.
         Comparing the explicit formula for $\bar g$ in \eqref{bar_tilde_g_def}
         with the notation \eqref{g_x_brackets} implies that
         $\bar L^*_{\ve\xi}g_{[d]}=\delta(\xi_1-d)\,\delta(\xi_2-x_2)\,\delta(\xi_3-x_3)$.
         So, by \eqref{eq:Green_adj}, $\bar L^*_{\ve\xi}g_{[x_1]}=L_{\ve\xi}^*G$,
         and also $\bar L^*_{\ve\xi}g_{[d]}=0$ for $d=-x_1,2\pm x_1$ and all $\ve\xi\in\Omega$.
         Now, by \eqref{bar_G_g}, we conclude that
         $\bar\phi=-\bar L^*_{\ve\xi}[\omega_1(\xi_1) \bar {\mathcal G}_2]$
         where $\bar {\mathcal G}_2:=\lambda^- g_{[2-x_1]}-p\,\lambda^{\!+} g_{[x_1+2]}$,
         and $\bar L^*_{\ve\xi}\bar {\mathcal G}_2=0$ for $\ve\xi\in\Omega$.

         From these observations,
         $\bar\phi=2\eps\omega_1'(\xi_1)\pt_{\xi_1} \bar {\mathcal G}_2
          +[\eps\omega_1''(\xi_1)-2q\omega_1'(\xi_1)]\bar {\mathcal G}_2$.
         The definition \eqref{cut_off} of $\omega_1$ implies that
         $\bar\phi$ vanishes at $\xi_1=0$ and for $\xi_1\ge\frac13$.
         This implies the desired assertion \eqref{eq:bar_phi_boundary}.
         Furthermore, we now get
         \begin{multline*}
             \norm{\bar \phi(\ve{x};\cdot)}{1,1\,;\Omega}
            +\norm{\pt_{x_2}\bar \phi(\ve{x};\cdot)}{1\,;\Omega}
            +\norm{\pt_{x_3}\bar \phi(\ve{x};\cdot)}{1\,;\Omega}\\
            \leq C \bigl(
                      \norm{\bar {\mathcal G}_2(\ve{x};\cdot)}{2,1\,;[0\frac13]\times\R}
                     +\norm{D_{x_2}\bar {\mathcal G}_2(\ve{x};\cdot)}{1,1\,;[0\frac13]\times\R}
                     +\norm{D_{x_3}\bar {\mathcal G}_2(\ve{x};\cdot)}{1,1\,;[0\frac13]\times\R}
                   \bigr).
         \end{multline*}
         Combining this with the bounds \eqref{lambda_g_1_3} for the terms $\lambda^\pm g_{[2\pm x_1]}$
         of $\bar {\mathcal G}_2$, and the observation that
         $|D_{x_2}p|+|D_{x_3}p|\le C|\pt_q p|\le C$ and $\pt_{\xi_k} p=0,\,k=1,2,3$, yields
         our assertions for $\bar\phi$ in \eqref{tilde_bar_w}.

    (ii) Now we prove the desired estimate \eqref{tilde_bar_w} for $\tilde \phi$.
         By \eqref{bar_tilde_G}, throughout this part of the proof we set $q=\frac12a(\ve\xi)\ge\frac12\alpha$.
         Comparing the notation \eqref{g_x_brackets} with
         the explicit formula for $\tilde g$ in \eqref{bar_tilde_g_def}, we rewrite
         \eqref{eq:Green_prim_const}
         as $\tilde L_{\ve{x}}g_{[x_1]}=\delta(\ve{x}-\ve\xi)$.
         So $\tilde L_{\ve{x}}g_{[x_1]}=L_{\ve{x}}G$, by~\eqref{eq:Green_prim}.
         Next, for each value $d=-x_1,2\pm x_1$ respectively set
         $s=-\xi_1,\mp(2-\xi_1)$. Now by \eqref{eq:def_g0}, one has
         $\widehat r_{[d]}=\sqrt{(s-x_1)^2+(\xi_2-x_2)^2+(\xi_3-x_3)^2}/\eps$
         so $g(\ve{x};s,\xi_2,\xi_3;q)=\ds\frac{1}{4\pi\eps^2}\frac{\E^{q(s-x_1)/\eps-q\widehat r_{[d]}}}{\widehat r_{[d]}}$.
         Note that $\tilde L_{\ve{x}}g(\ve{x};s,\xi_2,\xi_3;q)=\delta(x_1-s)\,\delta(x_2-\xi_2)\,\delta(x_3-\xi_3)$
         and none of our three values of $s$ is in $[0,1]$ (i.e. $\delta(s-x_1)=0$).
         Consequently,
         $\tilde L_{\ve{x}}
         \Bigl[\ds\frac{\E^{q(\widehat\xi_{1,[x_1]}-\widehat r_{[d]})}}
         {\widehat r_{[d]}}\Bigr]=0$
         for all $\ve{x}\in\Omega$.
         Comparing \eqref{tilde_G_def} and \eqref{tilde_G_g}, we now conclude that
         $\tilde\phi=-\tilde L_{\ve{x}}[\omega_0(\xi_1) \tilde {\mathcal G}_2]$
         where $\tilde{\mathcal G}_2:=p\, g_{[-x_1]}-p\,\lambda^{\!+} g_{[2+x_1]}$
         and $\tilde L_{\ve{x}}\tilde{\mathcal G}_2=0$ for $\ve{x}\in\Omega$.

         From these observations,
         $\tilde\phi =2\eps\omega_0'(x_1)\,\pt_{x_1} \tilde{\mathcal G}_2
         +[\eps\omega_0''(x_1)+2q\omega_0'(x_1)]\tilde{\mathcal G}_2$.
         As the definition~\eqref{cut_off} of $\omega_0$ implies that $\tilde\phi$
         vanishes for $x_1\le\frac23$, we have
         \[
            \norm{\tilde\phi(\ve{x};\cdot)}{1,1\,;\Omega}
            \leq C \max_{\stackrel{\ve{x}\in[\frac23,1]\times\R^2}
                                  {k\,=\,0,1}}
                   \norm{\pt^k_{x_1}\tilde {\mathcal G}_2(\ve{x};\cdot)}{1,1\,;\Omega}\,.
         \]
         Here $\tilde G_2$ is smooth and has no singularities for $x_1\in[\frac23,1]$
         (because $\widehat r_{[2+x_1]}\geq \widehat r_{[-x_1]}\geq \frac23\eps^{-1}$ for $x\in[\frac23,1]$).
         Note that $\norm{\pt^k_{x_1} g_{[-x_1]}}{1,1\,;\Omega}\le C \eps^{-2}$,
         and       $\norm{\pt^k_{x_1} (\lambda^{\!+}g_{[2+x_1]})}{1,1\,;\Omega}\le C \eps^{-2}$
         (these two estimates are similar to the ones in
         Lemmas~\ref{lem:g0_bounds} and~\ref{lem:g_lmbd_bounds},
         but easier to deduce as they are not sharp).
         We combine these two bounds with
         $|\pt^k_{x_1} \pt^l_{\xi_1}\pt^m_{\xi_2}\pt^n_{\xi_3} p|
          \leq C \eps^{-2}p
         = C \eps^{-2}\E^{-2qx_1/{\eps}}$ for $k,\, l+m+n\leq 1$.
         As for $x_1\ge\frac23$ we enjoy the bound
         $\E^{-2qx_1/{\eps}}
          \leq \E^{-\frac23\alpha/{\eps}}
          \leq  C\eps^4 \E^{-\frac12\alpha/{\eps}}$, the desired estimate for $\tilde\phi$ follows.
   \end{proof}

   \begin{lem}\label{lem:tilde_bar_G}
      Let the function $R=R(\ve{x};\ve\xi)$ be such that
      $|R|\le C\min\{\eps\widehat r_{[\ve{x}]},1\}$.
      The functions $\bar G$ and $\tilde G$ of \eqref{bar_tilde_G}, \eqref{bar_G_tilde_g}
      satisfy
      \begin{subequations}
         \begin{align}
            \norm{\bar G(\ve{x};\cdot)}{1\,;\Omega}+\norm{\tilde G(\ve{x};\cdot)}{1\,;\Omega}
               &\le C,\label{G_bar_t_1}\\
            \norm{\pt_{\xi_1}\bar G(\ve{x};\cdot)}{1\,;\Omega}
               &\le C(1+|\ln\eps|),\label{bar_G_xi}\\
            \norm{\pt_{\xi_k}\bar G(\ve{x};\cdot)}{1\,;\Omega}
               &\le C\eps^{-1/2},\quad k=2,\,3,\label{bar_G_eta}\\
            \norm{(R\,\pt_{\xi_1}\bar G)(\ve{x};\cdot)}{1\,;\Omega}
            +\eps^{1/2}\norm{(R\,\pt^2_{\xi_1\xi_k}\bar G)(\ve{x};\cdot)}{1\,;\Omega}
               &\le C,\quad\,\qquad k=2,\,3,\label{bar_G_xi_eta}
      \intertext{and for any ball $B(\ve{x}';\rho)$ of radius $\rho$ centered at any
               $\ve{x}'\in[0,1]\times\R^2$, one has}
            |\bar G(\ve{x};\cdot)|_{1,1\,;B(\ve{x}';\rho)\cap\Omega}
               &\le C\eps^{-1}\rho,\label{bar_G_1_1_B}
         \end{align}
         while for the ball $B(\ve{x};\rho)$ of radius $\rho$ centered at $\ve{x}$,
         we have
         \begin{align}
            \norm{\pt^2_{\xi_1}\bar G(\ve{x};\cdot)}{1\,;\Omega\setminus B(\ve{x};\rho)}
               &\le C\eps^{-1}\ln(2+\eps/\rho),\label{bar_G_xi_xi_2}\\
            \norm{\pt^2_{\xi_k}\bar G(\ve{x};\cdot)}{1\,;\Omega\setminus B(\ve{x};\rho)}
               &\le C\eps^{-1}(\ln(2+\eps/\rho)+|\ln\eps|),\quad k=2,\,3.\label{bar_G_eta_eta_2}
         \end{align}
         Furthermore, we have for $k=2,\,3$
         \begin{align}
             \norm{\pt_{x_k}\bar G(\ve{x};\cdot)}{1\,;\Omega}
            +\norm{(R\,\pt^2_{\xi_1 x_k}\bar G)(\ve{x};\cdot)}{1\,;\Omega}
               &\le C\eps^{-1/2},\label{bar_G_y}\\
            \norm{\pt_{\xi_k}\tilde G(\ve{x};\cdot)}{1\,;\Omega}
               &\le C\eps^{-1/2},\label{G_tilde_eta}\\
            \int_0^1\!\bigl(\norm{(R\,\pt^2_{x_1\xi_k}\tilde G)(\ve{x};\cdot)}{1\,;\Omega}+
                            \norm{\pt_{x_1}\tilde G(\ve{x};\cdot)}{1\,;\Omega}\bigr)\,dx_1
               &\le C\eps^{-1/2}.\label{G_tilde_eta_x}
         \end{align}
      \end{subequations}
   \end{lem}

   \begin{proof}
      Throughout the proof,
whenever $k$ appears in any relation,
it will be understood to be valid for $k=2,3$.

      First, note that $\widehat r_{[-x_1]}\geq \widehat r_{[x_1]}$
      and $\widehat r_{[2\pm x_1]}\ge \widehat r_{[x_1]}$ for all $\ve\xi\in\Omega$, therefore
      \begin{gather}\label{R_ext}
         |R|\leq C\,\min\bigl\{\eps\widehat r_{[x_1]},\,\eps\widehat r_{[-x_1]},\,\eps\widehat r_{[2-x_1]},\,\eps\widehat r_{[2+x_1]},
                               \,1\bigr\}.
      \end{gather}

      Note also that in view of Remark~\ref{rem_q_var},
      all bounds of Lemma~\ref{lem:g0_bounds}  apply to
      the components $g_{[\pm x_1]}$
      and all bounds of Lemma~\ref{lem:g_lmbd_bounds} apply to
      the components $\lambda^\pm g_{[2\pm x_1]}$
      of $\bar{\mathcal G}$ and $\tilde{\mathcal G}$ in \eqref{bar_G_tilde_g}.

\underline{\it Asterisk notation.}
      In some parts of this proof, when discussing derivatives of $\bar {\mathcal G}$,
      we shall use the  notation $\bar {\mathcal G}^\ast$ 
      prefixed by some differential operator, e.g., $\pt_{x_1} \bar {\mathcal G}^\ast$.
      This will mean that the differential operator is applied only to the terms of the
      type $g_{[d\pm x_1]}$, e.g., $\pt_{x_1} \bar {\mathcal G}^\ast$ is obtained by replacing
      each of the four terms $g_{[d\pm x_1]}$ in the definition \eqref{bar_G_g}
      of $\bar {\mathcal G}$ by $\pt_{x_1} g_{[d\pm x_1]}$ respectively.
      \begin{enumerate}
         \item The first desired estimate \eqref{G_bar_t_1}
               follows from the bound \eqref{g_L1} for $g_{[\pm x_1]}$
               and the bound \eqref{g_L1_mu} for $\lambda^\pm g_{[2\pm x_1]}$
               combined with $|p|\le 1$ and $|\omega_{0,1}|\le 1$
               (in fact, the bound for $\bar G$ can obtained
               by imitating the proof of Lemma~\ref{lem_G_L1}).
         \item Rewrite \eqref{bar_G_g} as
               \[
                  \bar {\mathcal G}=\bar {\mathcal G}_1-\omega_1(\xi_1)\bar {\mathcal G}_2,
                  \quad\mbox{where}\quad
                  \bar {\mathcal G}_1:=g_{[x_1]}-p\, g_{[-x_1]},
                  \quad
                  \bar {\mathcal G}_2:=\lambda^- g_{[2-x_1]}-p\,\lambda^{\!+} g_{[2+x_1]}.
               \]
               As $q=\frac12a(\ve{x})$ in $\bar G$
               (i.e. $p$ and $\lambda^\pm$ in $\bar G$ do not involve $\ve\xi$),
               one gets
               \[
                  \pt_{\xi_1} \bar G=\pt_{\xi_1}\bar {\mathcal G}^\ast-\omega_1'(\xi_1)\bar {\mathcal G}_2,
                  \qquad
                  \pt_{\xi_k} \bar G=\pt_{\xi_k}\bar {\mathcal G}^\ast,
                  \qquad
                  \pt^2_{\xi_1\xi_k} \bar G=\pt^2_{\xi_1\xi_k}\bar {\mathcal G}^\ast
                  -\omega_1'(\xi_1)\pt_{\xi_k}\bar {\mathcal G}^\ast_2.
               \]
               Now the desired estimate \eqref{bar_G_xi}
               follows from the bound \eqref{g_xi_L1} for $\pt_{\xi_1} g_{[\pm x_1]}$,
               the bound \eqref{g_L1_mu2} for $\lambda^\pm \,\pt_{\xi_1} g_{[2\pm x_1]}$,
               and the bound \eqref{g_L1_mu} for $\lambda^\pm g_{[2\pm x_1]}$.
               Similarly, our next assertion \eqref{bar_G_eta}
               follows from the bound \eqref{g_eta_L1} for $\pt_{\xi_k} g_{[\pm x_1]}$,
               and the bound \eqref{g_L1_mu3} for $\lambda^\pm \pt_{\xi_k}g_{[2\pm x_1]}$.

               The next estimate \eqref{bar_G_xi_eta} is deduced using
               \[
                  |R\,\pt_{\xi_1} \bar G|
                    \leq   |R \,\pt_{\xi_1}\bar {\mathcal G}_1^\ast|
                         +C|\pt_{\xi_1}\bar {\mathcal G}_2^\ast|
                         +C|\bar {\mathcal G}_2|,
                  \quad
                  |R\,\pt^2_{\xi_1\xi_k} \bar G|
                     \leq   |R\,\pt^2_{\xi_1\xi_k}\bar {\mathcal G}^\ast|
                          +C|\pt_{\xi_k}\bar {\mathcal G}^\ast_2|.
               \]
               Here, in view of \eqref{R_ext}, the term $R \,\pt_{\xi_1}\bar {\mathcal G}_1^\ast$
               is estimated using the bound \eqref{g_xi_R_L1} for
               $\eps\widehat r_{[\pm x_1]} \pt_{\xi_1} g_{[\pm x_1]}$,
               while the terms $R\,\pt^2_{\xi_1\xi_k}\bar {\mathcal G}^\ast$
               are estimated using the bound \eqref{g_eta_x_L1} for
               $\eps\widehat r_{[\pm x_1]}\pt^2_{\xi_1\xi_k}g_{[\pm x_1]}$
               and the bound \eqref{g_L1_mu3} for
               $\lambda^\pm \eps\widehat r_{[2\pm x_1]}\pt^2_{\xi_1\xi_k}g_{[2\pm x_1]}$.
%
               The remaining terms $\pt_{\xi_1}\bar {\mathcal G}_2^\ast$,
               $\bar {\mathcal G}_2$
               and $\pt_{\xi_k}\bar{\mathcal G}^\ast_2$
               appear in  $\pt_{\xi_1}\bar G$ and $\pt_{\xi_k}\bar G$, 
               so have been bounded when obtaining \eqref{bar_G_xi}, \eqref{bar_G_eta}.
         \item The next assertion \eqref{bar_G_1_1_B} is proved similarly to
               \eqref{bar_G_xi} and \eqref{bar_G_eta}, only using
               the bound \eqref{g_eta_ball} for $g_{[\pm x_1]}$
               and the bound \eqref{g_L1_mu5} for $\lambda^\pm g_{[2\pm x_1]}$.
         \item As $q=\frac12a(\ve{x})$ in $\bar G$, then
               $\pt^2_{\xi_m} \bar G=\pt^2_{\xi_m}\bar {\mathcal G}^\ast$, $m=1,\,2,\,3$,
               and the assertions \eqref{bar_G_xi_xi_2} and \eqref{bar_G_eta_eta_2}
               immediately follow from the bounds \eqref{g_xi2_L1} and \eqref{g_eta2_L1} for
               $\pt^2_{\xi_m}g_{[\pm x_1]}$ combined with
               the bounds \eqref{g_L1_mu6} for $\lambda^\pm \pt^2_{\xi_m}g_{[2\pm x_1]}$
               where $m=1,\,2,\,3$.
         \item As $q=\frac12a(\ve{x})$ in $\bar G$, so
               using the operator $D_{x_k}$ of \eqref{D_ops}, one gets 
               \begin{align*}
                  \pt_{x_k} \bar G
                   =D_{x_k}\bigl[g_{[x_1]}-p \,g_{[-x_1]}\bigr]^\ast
                    -\omega_1(\xi_1)\,
                    \bigl[D_{x_k}(\lambda^-  g_{[2-x_1]})-p\, D_{x_k}(\lambda^{\!+} g_{[2+x_1]})\bigr]\\
                  -{\ts\frac12} \pt_{x_k} a(\ve{x})\cdot \pt_q p\cdot
                        \bigl[  g_{[-x_1]}-\omega_1(\xi_1)\lambda^{\!+} g_{[2+x_1]}\bigr],
               \end{align*}
               where $|\pt_q p|\le C$ by \eqref{lmb_p_def}
               (and we used the previously defined asterisk notation).
               Now,  $\pt_{x_k} \bar G$ is estimated using the
               bound \eqref{D_g} for $D_{x_k}g_{[\pm x_1]}$ and
               the bound \eqref{D_g_lmbd} for $D_{x_k}(\lambda^\pm  g_{[2\pm x_1]})$.
               For the term $g_{[-x_1]}$ in $\pt_{x_k} \bar G$
               we use the bound \eqref{g_L1}, and for the term  $\lambda^{\!+} g_{[2+x_1]}$
               the bound \eqref{g_L1_mu}.
               Consequently, one gets the desired bound \eqref{bar_G_y}
               for $D_{x_k}\bar G^\ast$.

               To estimate $R\,\pt^2_{\xi_1 x_k} \bar G$, $k=2,\,3$, a calculation shows that
               \begin{align*}
                 \pt^2_{\xi_1 x_k} \bar G
                  =(D_{x_k}\pt_{\xi_1}) \bigl[g_{[x_1]}-p \,g_{[-x_1]}\bigr]^\ast
                  \!\!-\omega_1(\xi_1)\,
                  \bigl[D_{x_k}(\lambda^- \pt_{\xi_1} g_{[2-x_1]})-p\, D_{x_k}(\lambda^{\!+}\pt_{\xi_1} g_{[2+x_1]})\bigr]\\
                 -{\ts\frac12} \pt_{x_k} a(\ve{x})\cdot \pt_q p\cdot
                   \bigl[\pt_{\xi_1}  g_{[-x_1]}
                         -\omega_1(\xi_1)\lambda^{\!+}\pt_{\xi_1} g_{[2+x_1]}\bigr]
                  -\omega_1'(\xi_1)\,\pt_{x_k}\bar G_2,
               \end{align*}
               where $\bar G_2:=\bar{\mathcal G}_2\bigr|_{q=a(\ve{x})/2}$.
               The assertion \eqref{bar_G_y}  for $R\,\pt^2_{\xi_1 x_k} \bar G$
               is now deduced as follows. In view of \eqref{R_ext},
               we employ the bound \eqref{D_pt_g} for the terms
               $\eps\widehat r_{[\pm x_1]}D_{x_k}\pt_{\xi_1} g_{[\pm x_1]}$
               and the bound \eqref{D_pt_lmbd} for the terms
               $\eps \widehat r_{[2\pm x_1]}\,D_{x_k}(\lambda^\pm\, \pt_{\xi_1} g_{[2\pm x_1]})$.
               For the remaining terms (that appear in the second line)
               we use $|R|\le C$ and $|\pt_q p|\le C$.
               Then we combine the bound \eqref{g_xi_L1} for $\pt_{\xi_1}  g_{[-x_1]}$
               and the bound \eqref{g_L1_mu2} for $\lambda^{\!+}\pt_\xi g_{[2+x_1]}$.
               The term $\pt_{x_k}\bar G_2$ is a part of $\pt_{x_k}\bar G$, which was estimated above,
               so for $\pt_{x_k}\bar G_2$
               we have the same bound as for $\pt_{x_k}\bar G$ in \eqref{bar_G_y}.
               This observation completes the proof of the bound
               for $R\,\pt^2_{\xi_1 x_k} \bar G$ in \eqref{bar_G_y}.
         \item We now proceed to estimating derivatives of $\tilde G$,
               so $q=\frac12a(\ve\xi)$ in this part of the proof.
               Let $\tilde {\mathcal G}^\pm:=g_{[\pm x_1]}-\lambda^{\mp} g_{[2\mp x_1]}$.
               Then \eqref{tilde_G_g}, \eqref{lmb_p_def} imply that
               $\tilde {\mathcal G}=\tilde {\mathcal G}^+-p_0\tilde {\mathcal G}^-$, where
               $p_0:=\omega_0(x_1)\,p=\omega_0(x_1)\,\E^{-2qx_1/\eps}$.
               Note that
               \[
                  D_{\xi_k} p_0
                    ={\ts\frac12}\pt_{\xi_k} a(\ve\xi)\cdot
                     (-2x_1/\eps)\,p_0,
                  \qquad
                  \pt_{x_1} p_0=[\omega_0'(x_1)-(2q/\eps)\,\omega_0(x_1)]\,\E^{-2qx_1/{\eps}}.
               \]
               Combining this with $|(-2x_1/\eps)\,p_0|\le C \E^{-qx_1/{\eps}}$
               and $q\ge \frac12 \alpha$ yields
               \begin{gather}\label{p_0}
                  |D_{\xi_k} p_0|\le C,\qquad
                  \int_0^1\!\bigl(|\pt_{x_1} p_0|+|D_{\xi_k} \pt_{x_1} p_0|\bigr)\,dx_1
                  \leq \int_0^1\!\bigl(C\eps^{-1}\E^{-\frac12\alpha x_1/{\eps}}\bigr)\,dx_1
                  \leq C.
               \end{gather}
               Furthermore, we claim that
               \begin{gather}\label{tilde_cal_G}
                  \norm{\tilde {\mathcal G}^-}{1\,;\Omega}\le C,
                  \qquad
                  \norm{\pt_{x_1}\tilde {\mathcal G}^\pm}{1\,;\Omega}\le C(1+|\ln\eps|),
                  \qquad
                  \norm{D_{\xi_k}\tilde {\mathcal G}^\pm}{1\,;\Omega}\le C\eps^{-1/2}.
               \end{gather}
%
               Here the first estimate
               follows from  the bounds \eqref{g_L1}, \eqref{g_L1_mu} for
               the terms $g_{[-x_1]}$ and $\lambda^{\!+}g_{[2+x_1]}$.
               The estimate for $\pt_{x_1}\tilde {\mathcal G}^\pm$ in \eqref{tilde_cal_G} follows from
               the bound \eqref{g_x_L1} for $\pt_{x_1} g_{[\pm x_1]}$ and
               the bound \eqref{D_g_lmbd} for $\pt_{x_1}(\lambda^{\pm}g_{[2\pm x_1]})$.
               Similarly, the estimate for $D_{\xi_k}\tilde {\mathcal G}^\pm$
               in \eqref{tilde_cal_G} is obtained using
               the bound \eqref{D_g} for $D_{\xi_k} g_{[\pm x_1]}$ and
               the bound \eqref{D_g_lmbd} for $D_{\xi_k}(\lambda^{\pm}g_{[2\pm x_1]})$.

               Next, a calculation shows that 
               \[
                  \pt_{\xi_k}\tilde G
                    = D_{\xi_k}\tilde {\mathcal G}^+
                     -p_0\,D_{\xi_k}\tilde {\mathcal G}^-
                     -D_{\xi_k} p_0\cdot\tilde {\mathcal G}^-,
                  \qquad
                  \pt_{x_1}\tilde G
                    = \pt_{x_1}\tilde {\mathcal G}^+
                     -p_0\,\pt_{x_1}\tilde {\mathcal G}^-
                     -\pt_{x_1} p_0\cdot\tilde {\mathcal G}^-.
               \]
               Combining these with \eqref{p_0}, \eqref{tilde_cal_G} yields
               \eqref{G_tilde_eta} and the bound for $\pt_{x_1}\tilde G$ in \eqref{G_tilde_eta_x}.

               To establish the estimate for $R\,\pt^2_{x_1\xi_k}\tilde G$
               in \eqref{G_tilde_eta_x}, 
               note that
               \[
                  \pt^2_{x_1\xi_k}\tilde G
                    = D_{\xi_k}\pt_{x_1} \tilde {\mathcal G}^+
                     -p_0\cdot D_{\xi_k}\pt_{x_1} \tilde {\mathcal G}^-\!
                     -\pt_{x_1} p_0\cdot D_{\xi_k} \tilde {\mathcal G}^-\!
                     -\pt_{\xi_k} p_0\cdot \pt_{x_1} \tilde {\mathcal G}^-\!
                     -D_{\xi_k}\pt_{x_1} p_0\cdot  \tilde {\mathcal G}^-\!.
               \]
               In view of \eqref{R_ext}, \eqref{p_0} and \eqref{tilde_cal_G},
               it now suffices to show that
               $\|R\,D_{\xi_k}\pt_{x_1} \tilde {\mathcal G}^\pm\|_{1\,;\Omega}\le C\eps^{-1/2}$.
               This latter estimate immediately
               follows from the bound \eqref{D_pt_g}
               for the terms $\eps\widehat r_{[\pm x_1]}\, D_{\xi_k}\pt_{x_1} g_{[\pm x_1]}$
               and the bound \eqref{D_pt_lmbd} for the terms
               $\eps\widehat r_{[\pm x_1]}\, D_{\xi_k}\pt_{x_1} (\lambda^\pm g_{[2\pm x_1]})$.
               This completes the proof of \eqref{G_tilde_eta_x}.
      \end{enumerate}
   \end{proof}
%
   \subsection{Approximations for the Green's function $G$ in the domain $\Omega=(0,1)^3$}\label{ssec:bounds_green_const_2}
      We now define approximations, denoted by $\bar G_\cube$ and $\tilde G_\cube$,
      for the Green's function $G$ in
      our original domain $\Omega=(0,1)^3$. For this, we use
      the approximations $\bar G$ and $\tilde G$ of \eqref{bar_tilde_G}, \eqref{bar_tilde_G20}
      for the domain $(0,1)\times\R^2$
      and again employ the method of images with an inclusion of the cut-off
      functions of \eqref{cut_off} in a two-step process as follows:
      \begin{subequations}\label{eq:def_tilde_bar_G}
      \begin{align}
         \!\bar G_{\mbox{\tiny$\Box$}}(\ve{x};\ve\xi)
         &:=  \bar G(\ve{x};\ve\xi)
               \hspace{-0.0cm}&&-\omega_0(\xi_2)\,\bar G(\ve{x};\xi_1,-\xi_2,\xi_3)
               \hspace{0.1cm}&&-\omega_1(\xi_2)\,\bar G(\ve{x};\xi_1,2-\xi_2,\xi_3),\notag\\
         \!\bar G_\cube(\ve{x};\ve\xi)
         &:=  \bar G_{\mbox{\tiny$\Box$}}(\ve{x};\ve\xi)
               \hspace{-1cm}&&-\omega_0(\xi_3)\,\bar G_{\mbox{\tiny$\Box$}}(\ve{x};\xi_1,\xi_2,-\xi_3)
               \hspace{-1cm}&&-\omega_1(\xi_3)\,\bar G_{\mbox{\tiny$\Box$}}(\ve{x};\xi_1,\xi_2,2-\xi_3),\\
         \!\tilde G_{\mbox{\tiny$\Box$}}(\ve{x};\ve\xi)
         &:=  \tilde G(\ve{x};\ve\xi)
               \hspace{-1cm}&&-\omega_0(x_2)\,\tilde G(x_1,-x_2,x_3;\ve\xi)
               \hspace{-1cm}&&-\omega_1(x_2)\,\tilde G(x_1,2-x_2,x_3;\ve\xi),\notag\\
         \!\tilde G_\cube(\ve{x};\ve\xi)
         &:=  \tilde G_{\mbox{\tiny$\Box$}}(\ve{x};\ve\xi)
               \hspace{-1cm}&&-\omega_0(x_3)\,\tilde G_{\mbox{\tiny$\Box$}}(x_1,x_2,-x_3;\ve\xi)
               \hspace{-1cm}&&-\omega_1(x_3)\,\tilde G_{\mbox{\tiny$\Box$}}(x_1,x_2,2-x_3;\ve\xi).
      \end{align}
      \end{subequations}
      Then $\bar G_\cube\bigr|_{\xi_1=0,1}=0$
      and $\tilde G_\cube\bigr|_{x_1=0,1}=0$
      (as this is valid for $\bar G$ and $\tilde G$, respectively),
      and furthermore, by \eqref{cut_off}, we have
      $\bar G_\cube\bigr|_{\xi_k=0,1}=0$
      and $\tilde G_\cube\bigr|_{x_k=0,1}=0$
      for $k=2,3$.
      \begin{rem}\label{rem:strip_to_square}
         Lemmas~\ref{lem_tilde_bar_w} and \ref{lem:tilde_bar_G}
         of the previous section remain valid if $\Omega$ is understood as $(0,1)^3$, and
         $\bar G$ and $\tilde G$ are replaced by $\bar G_\cube$ and $\tilde G_\cube$,
         respectively, in the definition \eqref{tilde_w_def} of $\bar\phi$ and $\tilde\phi$
         and in the lemma statements.

         This is shown by imitating the proofs of these two lemmas.
         We leave out the details and only note that the application of the method of images
         in the $\xi_2$- and $\xi_3$- ($x_2$- and $x_3$-) directions
         is relatively straightforward as an inspection of \eqref{eq:def_g0}
         shows that in these directions, the fundamental solution $g$
         is symmetric and exponentially decaying away from the singular point.
      \end{rem}
      As $\bar G_\cube$ and $\tilde G_\cube$ in the domain $\Omega=(0,1)^3$
      enjoy the same properties as $\bar G$ and $\tilde G$ in the domain
      $(0,1)\times\R^2$, we shall sometimes skip the subscript
      \cube{} when there is no ambiguity.
%
\section{Proof of Theorem~\ref{thm:main} for $\Omega=(0,1)^3$\\
        (general variable-coefficient case)}\label{sec:main_proof}
%
We are now ready to establish our main result, Theorem~\ref{thm:main},
for the original variable-coefficient problem \eqref{eq:Lu} in the domain
$\Omega=(0,1)^3$.
In Section~\ref{sec:bounds_green_const}, we have already
obtained various bounds for the approximations $\tilde G_\cube$
and $\bar G_\cube$ of $G$ in $\Omega=(0,1)^3$.
So now we consider the two functions
\[
   \tilde v(\ve{x};\ve\xi):=[G-\tilde G_\cube](\ve{x};\ve\xi),
   \qquad
   \bar v(\ve{x};\ve\xi)=[G-\bar G_\cube](\ve{x};\ve\xi).
\]
Throughout this section, we shall skip the subscript \cube{}
as we always deal with the domain $\Omega=(0,1)^3$.

Note that, by \eqref{tilde_w_def}, we have
$L_{\ve{x}}\tilde v=L_{\ve{x}}[G-\tilde G]
   = [\tilde L_{\ve{x}}-L_{\ve{x}}]\tilde G-\tilde\phi$,
and similarly
$L^*_{\ve\xi}\bar v=L^*_{\ve\xi}[G-\bar G]
   = [\bar L^*_{\ve\xi}-L^*_{\ve\xi}]\bar G-\bar\phi$.
Consequently, the functions $\tilde v$ and $\bar v$ are solutions of the following problems:
\begin{subequations}
   \begin{align}
      L_{\ve{x}} \tilde v(\ve{x};\ve\xi)
         &= \tilde h(\ve{x};\ve\xi)\;\;\mbox{for}\;\ve{x}\in\Omega,
            \qquad \tilde v(\ve{x};\ve\xi)=0\;\;\mbox{for}\; \ve{x}\in\pt\Omega,\label{prob_tilde_v}\\
      L^*_{\ve\xi} \bar v(\ve{x};\ve\xi)
         &= \bar h(\ve{x};\ve\xi)\;\;\mbox{for}\;\ve\xi\in\Omega,
            \qquad \bar v(\ve{x};\ve\xi)=0\;\;\mbox{for}\; \ve\xi\in\pt\Omega.\label{prob_bar_v}
   \end{align}
\end{subequations}
Here the right-hand sides are given by
\begin{subequations}
   \begin{align}
      \tilde h(\ve{x};\ve\xi)
         &:= \pt_{x_1}\{R\,\tilde G\}(\ve{x};\ve\xi)
            -b(\ve{x})\,\tilde G(\ve{x};\ve\xi)
            -\tilde \phi(\ve{x};\ve\xi),\label{tilde_bar_h}\\
      \bar h(\ve{x};\ve\xi)
         &:= \{R\,\pt_{\xi_1}\bar G\}(\ve{x};\ve\xi)
            -b(\ve\xi)\,\bar G(\ve{x};\ve\xi)
            -\bar \phi(\ve{x};\ve\xi),\label{bar_h}
   \end{align}
\end{subequations}
where
\begin{gather}\label{R_def}
   R(\ve{x};\ve\xi):=a(\ve{x})-a(\ve\xi),\qquad\mbox{so}\;\;
   |R|\le C\min\{\eps\widehat r_{[x_1]},1\}.
\end{gather}
Applying the solution representation formulas \eqref{eq:sol_prim} and
\eqref{eq:sol_adj} to problems \eqref{prob_tilde_v} and \eqref{prob_bar_v},
respectively, one gets
\begin{subequations}
   \begin{align}
      \tilde v(\ve{x};\ve\xi)
         &= \iiint_\Omega G(\ve{x};\ve{s})\,\tilde h(\ve{s};\ve\xi)\,d\ve{s},\label{tilde_v}\\
      \bar v(\ve{x};\ve\xi)
         &= \iiint_\Omega G(\ve{s};\ve\xi)\,\bar h(\ve{x};\ve{s})\,d\ve{s}.\label{bar_v}
   \end{align}
\end{subequations}
We now proceed to the {proof of Theorem~\ref{thm:main}}.
\smallskip

\begin{proof}
       Throughout the proof,
whenever $k$ appears in any relation,
it will be understood to be valid for $k=2,3$.

(i)
   First we establish \eqref{eq:thm:G_eta}. 
   Note that, the bounds \eqref{G_tilde_eta} and \eqref{bar_G_y}
   for $\pt_{\xi_k}\tilde G$ and $\pt_{x_k}\bar G$,
   respectively, it suffices to show that
   $ \norm{\pt_{\xi_k}\tilde v(\ve{x};\cdot)}{1\,;\Omega}
    +\norm{\pt_{x_k}\bar v(\ve{x};\cdot)}{1\,;\Omega}\le C\eps^{-1/2}$.

   Applying $\pt_{\xi_k}$ to \eqref{tilde_v} and $\pt_{x_k}$ to \eqref{bar_v}, we arrive at
   \begin{align*}
     \pt_{\xi_k}\tilde v(\ve{x};\ve\xi)
        &= \iiint_\Omega \!G(\ve{x};\ve{s})\,\pt_{\xi_k}\tilde h(\ve{s};\ve\xi)\,d\ve{s},\\
     \pt_{x_k}\bar v(\ve{x};\ve\xi)
        &=\iiint_\Omega\!G(\ve{s};\ve\xi)\,\pt_{x_k}\bar h(\ve{x};\ve{s})\,d\ve{s}.
   \end{align*}
   From this, a calculation shows that
   \begin{align*}
      \norm{\pt_{\xi_k}\tilde v(\ve{x};\cdot)}{1\,;\Omega}
         &\le \Bigl(\sup_{s_1\in(0,1)}\iint_{\R^2} |G(\ve{x};\ve{s})|\,ds_2\,ds_3\Bigr)\cdot
               \int_0^1\!\!\sup_{(s_2,s_3)\in\R^2}
               \norm{\pt_{\xi_k} \tilde h(\ve{s};\cdot)}{1\,;\Omega}\,ds_1\,,\\
      \norm{\pt_{x_k}\bar v(\ve{x};\cdot)}{1\,;\Omega}
         &\le \Bigl(\,\sup_{\ve{s}\in\Omega}\norm{G(\ve{s};\cdot)}{1\,;\Omega}.\Bigr)\cdot
               \norm{\pt_{x_k}\bar h(\ve{x};\cdot)}{1\,;\Omega}.
   \end{align*}
   So, in view of \eqref{G_L1}, to prove \eqref{eq:thm:G_eta}, it remains to show that
   \[
      \int_0^1\!\!\sup_{(x_2,x_3)\in\R^2}\norm{\pt_{\xi_k} \tilde h(\ve{x};\cdot)}{1\,;\Omega}\,dx_1
        \leq C\eps^{-1/2},
      \qquad
      \norm{\pt_{x_k}\bar h(\ve{x};\cdot)}{1\,;\Omega}
         \le C\eps^{-1/2}.
   \]
   These two bounds follow from the definitions \eqref{tilde_bar_h}, \eqref{R_def}
   of $\tilde h$ and $\bar h$, which imply that
   \begin{align*}
      |\pt_{\xi_k} \tilde h(\ve{x};\ve\xi)|
         &\le |R\,\partial^2_{x_1 \xi_k}\tilde G|
             +C\bigl( |\pt_{x_1}\tilde G|+|\pt_{\xi_k}\tilde G|\bigr)
             +|\pt_{\xi_k} \tilde\phi|,\\
      |\pt_{x_k}\bar h(\ve{x};\ve\xi)|
         &\le |R\,\pt^2_{\xi_1 \xi_k}\bar G|
             +C\bigl(|\pt_{\xi_1}\bar G|+ |\pt_{x_k}\bar G|\bigr)
             +|\pt_{x_k}\bar\phi|,
   \end{align*}
   combined with the bounds \eqref{tilde_bar_w} for $\bar\phi$, $\tilde\phi$,
   the bounds \eqref{G_tilde_eta}, \eqref{G_tilde_eta_x} for $\tilde G$
   and the bounds \eqref{bar_G_xi}, \eqref{bar_G_y} for $\bar G$.
   Thus we have shown \eqref{eq:thm:G_eta}. 
   \smallskip

(ii)
   Next we proceed to obtaining
   the assertions \eqref{eq:thm:G_xi},
   \eqref{eq:thm:G_xixi}  and \eqref{eq:thm:G_etaeta}.
   We claim that to get these bounds, it suffices to show that
   \begin{subequations}\label{desired}
   \begin{align}
     \mathcal{V}
      :=\max_{k=2,\,3}\,\,\sup_{\ve{x}\in\Omega}\norm{\pt^2_{\xi_k}\bar v(\ve{x};\cdot)}{1\,;\Omega}
      &\le C(\eps^{-1}+\eps^{-1/2}\, \mathcal{W}),
      \label{bar_v_eta_eta}\\
     \mathcal{W}
      :=\sup_{\ve{x}\in\Omega}\norm{\partial_{\xi_1} G(\ve{x};\cdot)}{1\,;\Omega}
      &\le C(1+|\ln\eps|+\eps\mathcal{V}),\label{mathcal_G}\\
      \sup_{\ve{x}\in\Omega}\|\pt^2_{\xi_1}\bar v(\ve{x};\cdot)\|_{1\,;\Omega}
      &\le C\,\eps^{-1}(1+\eps\mathcal{V}).\label{bar_v_xi_xi}
   \end{align}
   \end{subequations}
   Indeed, there is a sufficiently small constant $c_*$ such that
   for $\eps\leq c_*$, combining the bounds \eqref{bar_v_eta_eta},\,\eqref{mathcal_G},
   one gets $\mathcal{W}\le C(1+|\ln\eps|)$, which is identical with
   \eqref{eq:thm:G_xi}.
   Then \eqref{bar_v_eta_eta} implies that $\mathcal{V}\le C\eps^{-1}$,
   which, combined with \eqref{bar_G_eta_eta_2},
   yields \eqref{eq:thm:G_etaeta}.
   Finally, $\mathcal{V}\le C\eps^{-1}$ combined with
   \eqref{bar_v_xi_xi} and then
   \eqref{bar_G_xi_xi_2} yields \eqref{eq:thm:G_xixi}.

   In the simpler non-singularly-perturbed case of $\eps>c_*$,
   by imitating part (i) of this proof, one obtains $\mathcal{W}\le C_1$,
   where $C_1$ depends on $c_*$.
   Combining this bound with \eqref{bar_v_eta_eta} and \eqref{bar_v_xi_xi},
   we again get \eqref{eq:thm:G_xi}, \eqref{eq:thm:G_xixi} and
   \eqref{eq:thm:G_etaeta}.

   We shall obtain \eqref{bar_v_eta_eta} in part (iii)
   and \eqref{mathcal_G} with \eqref{bar_v_xi_xi} in part (iv) below.
   \smallskip

(iii)
   To get \eqref{bar_v_eta_eta},
    it suffices to set $k=2$ and
   consider $\bar V:=\pt^2_{\xi_2}\bar v$ (as $\pt^2_{\xi_3}\bar v$
   is estimated similarly).
   The problem \eqref{prob_bar_v} for $\bar v$ implies that
   \begin{gather}\label{bar_V}
      L^*_{\ve\xi}\bar V(\ve{x};\ve\xi)=\bar H(\ve{x};\ve\xi)
         \;\;\mbox{for}\;\ve\xi\in\Omega,\quad
      \bar V(\ve{x};\ve\xi)=0\;\;\mbox{for}\;\ve\xi\in\partial\Omega.
   \end{gather}
    The homogeneous boundary conditions $\pt^2_{\xi_2}\bar v\bigr|_{\xi_m=0,1}=0$
   in \eqref{bar_V}
for $m=1,\,3$
   immediately follow from $\bar v\bigr|_{\xi_m=0,1}=0$.
   The homogeneous boundary conditions on the boundary edges
   $\xi_2=0,1$ are obtained as follows. As $\bar v\bigr|_{\xi_2=0,1}=0$ so
   $\pt_{\xi_1}\bar v\bigr|_{\xi_2=0,1}=\pt^2_{\xi_m}\bar v\bigr|_{\xi_2=0,1}=0$,
   where again $m=1,\,3$.
   Combining this with $\bar h\bigr|_{\xi_2=0,1}=0$
   (for which, in view of Remark~\ref{rem:strip_to_square}, we used \eqref{eq:bar_phi_boundary})
   and the differential equation for $\bar v$ at $\xi_2=0,1$, one
   finally gets
   $\pt^2_{\xi_2}\bar v\bigr|_{\xi_2=0,1}=0$.

   For the right-hand side $\bar H$ in \eqref{bar_V}, a calculation shows
   that $\bar H=\bar H(\ve{x};\ve\xi)=\pt_{\xi_2}\bar h_1+\bar h_2$
   with $\bar h_{1}=\bar h_{1}(\ve{x};\ve\xi)$ and
   $\bar h_{2}=\bar h_{2}(\ve{x};\ve\xi)$ defined by
   \[
      \bar h_1:=\partial_{\xi_k}\bar h
                -2\pt_{\xi_k} a(\ve\xi)\cdot\pt_{\xi_1}\bar v,
      \qquad
      \bar h_2:=\pt^2_{\xi_k} a(\ve\xi)\cdot\pt_{\xi_1}\bar v
                -2\pt_{\xi_k} b(\ve\xi)\cdot\pt_{\xi_k}\bar v
                -\pt^2_{\xi_k} b(\ve\xi)\cdot\bar v,
   \]
with $k=2$.
   Here we used
   $ \pt^2_{\xi_k}[a\,\pt_{\xi_1}\bar v]
    =a\,\pt_{\xi_1}\! \bar V+2\pt_{\xi_k} a \,\pt^2_{\xi_1\xi_k}\bar v
     +\pt^2_{\xi_k} a\,\pt_{\xi_1}\bar v
    =a\,\pt_{\xi_1}\! \bar V+\pt_{\xi_k}[2\,\pt_{\xi_k} a \,\pt_{\xi_1}\bar v]
     -\pt^2_{\xi_k} a\,\pt_{\xi_1}\bar v$
   and
   $\pt^2_{\xi_k}[b\bar v]=b \bar V+2\,\pt_{\xi_k} b \,\pt_{\xi_k}\bar v+\pt^2_{\xi_k} b\,\bar v$.
   (Note that  $\bar H$ is understood in the sense of distributions; see Remark \ref{rem_H} below.)

   Now, applying the solution representation formula \eqref{eq:sol_adj} to problem \eqref{bar_V},
   and then integrating the term with $\bar h_1$ by parts, yields
   \[
      \bar V(\ve{x};\ve\xi)
       =\iiint_\Omega\!\bigl[-\pt_{s_2} G(\ve{s};\ve\xi)\,\bar h_1(\ve{x};\ve{s})
        +
        G(\ve{s};\ve\xi)\,\bar h_2(\ve{x};\ve{s})\bigr]\,d\ve{s},
   \]
   (for the validity of the above integration by parts we again refer
   to Remark~\ref{rem_H}).
   As \eqref{eq:thm:G_eta} implies $\sup_{\ve{s}\in\Omega}\norm{\pt_{s_2} G(\ve{s};\cdot)}{1\,;\Omega}\le C\eps^{-1/2}$,
   while \eqref{G_L1} implies $\sup_{\ve{s}\in\Omega}\|G(\ve{s};\cdot)\|\le C$,
   imitating the argument used in part (i) of this proof yields
   \[
      \norm{\pt^2_{\xi_2}\bar v(\ve{x};\cdot)}{1\,;\Omega}
      =   \norm{\bar V(\ve{x};\cdot)}{1\,;\Omega}
      \leq C\bigl( \eps^{-1/2}\norm{\bar h_1(\ve{x};\cdot)}{1\,;\Omega}
                 +\norm{\bar h_2(\ve{x};\cdot)}{1\,;\Omega}\bigr).
   \]
   So to get our assertion \eqref{bar_v_eta_eta}, it remains to show that
   \begin{equation}\label{h12_bound}
     \norm{\bar h_{1}(\ve{x};\cdot)}{1\,;\Omega}+\norm{\bar h_{2}(\ve{x};\cdot)}{1\,;\Omega}\le C(\eps^{-1/2}+\mathcal{W}).
   \end{equation}
   To check this latter bound, note that
   $|\bar h_1|+|\bar h_2|
    \le C(|\pt_{\xi_k}\bar h|+|\pt_{\xi_1}\bar v|+|\pt_{\xi_k}\bar v|+|\bar v|)$
    with $k=2$.
   Note also that
   \[
      \norm{\bar v(\ve{x};\cdot)}{1,1\,;\Omega}
         \le C(\eps^{-1/2}+\mathcal{W})+\norm{\bar G(\ve{x};\cdot)}{1,1\,;\Omega},
   \]
   where we employed $\bar v= G-\bar G$ and then the bounds \eqref{G_L1},
   \eqref{eq:thm:G_eta} and the definition \eqref{mathcal_G} of $\mathcal{W}$ for $G$.
   Combining these two observations with
   \[
      |\pt_{\xi_k}\bar h(\ve{x};\ve\xi)|
       \le |R\,\pt^2_{\xi_1\xi_k} \bar G|
          +C\bigl(|\pt_{\xi_1}\bar G|+|\pt_{\xi_k}\bar G|+|\bar G|\bigr)
          +|\pt_{\xi_k}\bar \phi|,
                   \qquad k=2,
   \]
   (where we used \eqref{bar_h}, \eqref{R_def}), and then with the bounds
   \eqref{G_bar_t_1}--\eqref{bar_G_xi_eta} for $\bar G$, and the bound
   \eqref{tilde_bar_w} for $\bar\phi$, one gets the required estimate~\eqref{h12_bound}.
   Thus \eqref{bar_v_eta_eta} is established.
   \smallskip

   (iv)
   To prove \eqref{mathcal_G} and \eqref{bar_v_xi_xi}, rewrite the problem \eqref{prob_bar_v} as a two-point
   boundary-value problem in $\xi_1$, in which $\ve{x}$, $\xi_2$ and $\xi_3$ appear
   as parameters, as follows
   \begin{gather}\label{eq:ode_barv}
      [-\eps\pt^2_{\xi_1}+a(\ve\xi)\pt_{\xi_1}]\,\bar v(\ve{x};\ve\xi)
      =\bar{\bar h}(\ve{x};\ve\xi)
      \quad\mbox{for}\;\xi_1\in(0,1),
      \qquad \bar v(\ve{x};\ve\xi)\bigr|_{\xi_1=0,1}=0,
   \end{gather}
   where
   \begin{gather}\label{bar_bar_h}
      \bar{\bar h}(\ve{x};\ve\xi)
        :=\bar h(\ve{x};\ve\xi)
           +\eps\,\bigl[\pt^2_{\xi_2}\bar v(\ve{x};\ve\xi)
           +\pt^2_{\xi_3}\bar v(\ve{x};\ve\xi)\bigr]
           -b(\ve\xi)\,\bar v(\ve{x};\ve\xi).
   \end{gather}
   Consequently, one can represent $\bar v$ via the Green's function
   $\Gamma=\Gamma(\xi_1,\xi_2,\xi_3;s)$
   of the one-dimensional operator $[-\eps\pt^2_{\xi_1}+a(\ve\xi)\pt_{\xi_1}]$.
   Note that $\Gamma$, for any fixed $\xi_2$, $\xi_3$ and $s$, satisfies the equation
   $[-\eps\pt^2_{\xi_1}+a(\ve\xi)\pt_{\xi_1}]\Gamma(\ve\xi;s)=\delta(\xi_1-s)$
   and the boundary conditions $\Gamma(\ve\xi;s)\bigr|_{\xi_1=0,1}=0$.
   Note also that
   \begin{gather}\label{G_1d}
      \int_{0}^{1}\!\! |\pt_{\xi_1} \Gamma(\ve\xi;s)|\,d\xi_1
       \le 2\alpha^{-1}
   \end{gather}
   \cite[Lemma 2.3]{And_MM02}; see also
   \cite[(I.1.18)]{RST08}, \cite[(3.10b) and Section~3.4.1.1]{Linss10}.

   The solution representation for $\bar v$ via $\Gamma$ is given by
   \[
      \bar v(\ve{x};\ve\xi)=\int_0^1\! \Gamma(\ve\xi;s)\, \bar{\bar h}(\ve{x};s,\xi_2,\xi_3)\,ds.
   \]
   Applying $\pt_{\xi_1}$ to this representation yields
   \[
      \norm{\pt_{\xi_1} \bar v(\ve{x};\cdot)}{1\,;\Omega}
      \le\left(\sup_{(s,\xi_2,\xi_3)\in \Omega}\int_{0}^{1}\!\! |\pt_{\xi_1} \Gamma(\ve\xi;s)|d\xi_1\right)
         \cdot\bignorm{\bar{\bar h}(\ve{x};\cdot)}{1\,;\Omega}.
   \]
   In view of \eqref{G_1d}, we now have
   $\norm{\pt_{\xi_1} \bar v}{1\,;\Omega} \le 2 \alpha^{-1} \norm{\bar{\bar h}}{1\,;\Omega}$.
   Note that the differential equation \eqref{eq:ode_barv} for $\bar v$
   implies that
   $\eps\norm{\pt_{\xi_1}^2\bar v}{1;\Omega}
     \leq C(\norm{\pt_{\xi_1}\bar v}{1;\Omega}+
                    \norm{\bar{\bar h}}{1;\Omega})$.
   So, furthermore, we get
   \[
      \| \pt_{\xi_1} \bar v\|_{1\,;\Omega}
      +\eps\| \pt^2_{\xi_1} \bar v\|_{1\,;\Omega} \le C \|\bar{\bar h}\|_{1\,;\Omega}.
   \]
   As $G=\bar v+\bar G$ and we have the bound \eqref{bar_G_xi} for $\pt_{\xi_1} \bar G$,
   to obtain the desired bounds \eqref{mathcal_G} and \eqref{bar_v_xi_xi}, it remains to show that
   $\norm{\bar{\bar h}(\ve{x};\cdot)}{1\,;\Omega}\le C+\eps\mathcal{V}$.
   Furthermore, the definitions \eqref{bar_bar_h} of $\bar{\bar h}$
   and \eqref{bar_v_eta_eta} of $\mathcal{V}$,
   imply that it now suffices to prove the two estimates
   \begin{gather}\label{two_bounds}
      \norm{\bar v(\ve{x};\cdot)}{1\,;\Omega}\le C,
      \qquad
      \norm{{\bar h}(\ve{x};\cdot)}{1\,;\Omega}\le C.
   \end{gather}
   The first of them follows from $\bar v=G-\bar G$
   combined with \eqref{G_L1} and \eqref{G_bar_t_1}.
   The second is obtained from the definition \eqref{bar_h} of ${\bar h}$
   using \eqref{bar_G_y} for $\norm{R\pt_{\xi_1}\bar G}{1\,;\Omega}$,
   \eqref{G_bar_t_1} for $\norm{\bar G}{1\,;\Omega}$
   and \eqref{tilde_bar_w} for $\norm{\bar \phi}{1\,;\Omega}$.
   This completes the proof of \eqref{mathcal_G} and \eqref{bar_v_xi_xi}, and thus of
   \eqref{eq:thm:G_xi}, \eqref{eq:thm:G_xixi} and \eqref{eq:thm:G_etaeta}.
   \smallskip

(v)
   We now focus on the remaining
%
   assertion \eqref{eq:thm:G_grad},
   again rewrite the problem \eqref{prob_bar_v} as
   \[
      [-\eps\laplace_{\ve\xi}+1]\,\bar v(\ve{x};\ve\xi)
      ={\bar h}_0(\ve{x};\ve\xi)\quad\mbox{for}\;\ve\xi\in\Omega,
      \qquad \bar v(\ve{x};\ve\xi)\bigr|_{\pt\Omega}=0,
   \]
   where
   \begin{gather}\label{bar_bar_h_new}
      {\bar h}_0(\ve{x};\ve\xi)
       :=  \bar h(\ve{x};\ve\xi)
          -a(\ve\xi)\,\pt_{\xi_1} \bar v(\ve{x};\ve\xi)
          +[1-b(\ve\xi)]\,\bar v(\ve{x};\ve\xi).
   \end{gather}
   We shall represent $\bar v$ via the Green's function $\Psi$
   of the two-dimensional self-adjoint operator
   $[-\eps\laplace_{\ve\xi}+1]$.
   Note that $\Psi=\Psi(\ve{s};\ve\xi)$, for any fixed $\ve{s}$, satisfies the equation
   $[-\eps\laplace_{\ve\xi}+1]\Psi(\ve{s};\ve\xi)=\delta(\ve\xi-\ve{s})$,
   and also the boundary conditions $\Psi(\ve{s};\ve\xi)\bigr|_{\ve\xi\in\pt\Omega}=0$.
   Furthermore, for any ball $B(\ve{x}';\rho)$ of radius $\rho$ centred at any
   $\ve{x}'$, we cite the estimate \cite[(3.5b)]{CK09}
   \begin{gather}\label{G_2d}
      |\Psi(\ve{s};\cdot)|_{1,1\,;B(\ve{x}';\rho)\cap\Omega} \le C\eps^{-1}\rho.
   \end{gather}
   The solution representation for $\bar v$ via $\Psi$ is given by
   \[
      \bar v(\ve{x};\ve\xi)
        =\iiint_\Omega \Psi(\ve{s};\ve\xi)\,{\bar h}_0(\ve{x};\ve{s})\,d\ve{s}.
   \]
   Applying $\pt_{\xi_m}$, $m=1,\,2,\,3$ to this representation yields
   \begin{gather}\label{ber_v_new}
      |\bar v(\ve{x};\cdot)|_{1,1\,;B(\ve{x}';\rho)\cap\Omega}
      \le \Bigl(\sup_{\ve{s}\in\Omega}|\Psi(\ve{s};\cdot)|_{1,1\,;B(\ve{x}';\rho)\cap\Omega}\Bigr)\cdot
          \norm{{\bar h}_0(\ve{x};\cdot)}{1\,;\Omega}.
   \end{gather}
   To estimate $\norm{{\bar h}_0}{1\,;\Omega}$, recall that
   it was shown in part (iv) of this proof that
   $\norm{\pt_{\xi_1} \bar v}{1\,;\Omega} \le 2 \alpha^{-1} \norm{\bar{\bar h}}{1\,;\Omega}$
   and
   $\norm{\bar{\bar h}(\ve{x};\cdot)}{1\,;\Omega}\le C+\eps\mathcal{V}$,
   and in part (ii) that $\mathcal{V}\le C\eps^{-1}$.
   Consequently $\norm{\pt_{\xi_1} \bar v}{1\,;\Omega} \le C$.
   Combining this with \eqref{bar_bar_h_new} and \eqref{two_bounds}
   yields $\norm{{\bar h}_0}{1\,;\Omega}\le C$.
   In view of \eqref{ber_v_new} and \eqref{G_2d},
   we now get $|\bar v|_{1,1\,;B(\ve{x}';\rho)\cap\Omega}\le C\eps^{-1}\rho$,
   which, combined with \eqref{bar_G_1_1_B},
   immediately gives the final desired bound \eqref{eq:thm:G_grad}.
 \end{proof}

\begin{rem}\label{rem_H}
   Note that the term $\pt^2_{\xi_k}\bar h$ in $\bar H$, where $k=2,\,3$,
   has such a singularity at $\ve\xi=\ve{x}$
   that it is not absolutely integrable in $\Omega$. So $\bar H$ and
   the differential equation in \eqref{bar_V} are understood in the sense of distributions
   \cite[Chapters 1,\,3]{Griffel}. In particular $\pt^2_{\xi_k}\bar h$ is a generalised
   $\xi_k$-derivative of the regular function $\pt_{\xi_k}\bar h$.
\end{rem}

%
   \bibliographystyle{plain}
   \bibliography{cd3d}
\end{document}